%% file: ramanujan.tex
\documentclass[12pt]{amsart}
\usepackage{amsmath,amsfonts,amssymb,amscd}
\usepackage{latexsym}

\usepackage[
hyperindex=true,pagebackref=true,bookmarks=true,
colorlinks=true,linkcolor=blue,citecolor=red]
{hyperref}
\def\~{{\rm --}} 

\font\dfont=cmbx10 at 11pt   

\title [Rogers-Ramanujan type identities and Nil-DAHA]
{Rogers-Ramanujan type identities and Nil-DAHA}
\author[Ivan Cherednik]{Ivan Cherednik $^\dag$}
\author[Boris Feigin]{Boris Feigin}
\date{October 28, 2012}

\thanks{$^\dag$  \date\ \ \ Partially supported by NSF grant
DMS--1101535}

\address[I. Cherednik]{Department of Mathematics, UNC
Chapel Hill, North Carolina 27599, USA\\
chered@email.unc.edu}
\address[B. Feigin]
{National Research University, Higher School of Economics, Russia,
Moscow, 101000, Myasnitskaya ul., 20;\\
Landau Institute for Theoretical Physics, Russia, Chernogolovka,
142432, prosp. Akademika Semenova, 1a;\\
Independent University of Moscow, Russia, Moscow, 119002, Bolshoy
Vlasyevskiy per., 11\\
borfeigin@gmail.com}

 \def\bysame{{\bf --- }}
 \def\~{{\bf --}}
\newcommand{\comment}[1]{}
\renewcommand{\Join}{\hbox{\small t\!o\!t}}
\renewcommand{\tilde}{\widetilde}
\renewcommand{\hat}{\widehat}

\input macro

\begin{document}
\maketitle
{\em\small Key words: Rogers-Ramanujan identities; Hecke algebras; 
$q$-Hermite polynomials; dilogarithm; Kac-Moody algebras;
Demazure characters; modular functions; coset algebras.}
\renewcommand{\baselinestretch}{1.2} 
{\textmd
\tableofcontents
} 
\renewcommand{\baselinestretch}{1.0}.
\vfill\eject

\renewcommand{\natural}{\wr}

\setcounter{section}{-1}
\setcounter{equation}{0}
\section{Introduction}
\comment{
In the theory of the Nil-DAHA Fourier transform, the inner products 
of q-Hermite polynomials for the measure function multiplied by a 
level-one theta function are the key. They are used to obtain 
expansions of products of any number of such theta functions in 
terms of the q-Hermite polynomials. An ample family of modular 
functions satisfying Rogers-Ramanujan type identities for arbitrary 
(reduced, twisted) affine root systems is obtained as an application. 
A relation to Rogers dilogarithm and Nahm's conjecture is discussed. 
Some of our q-series can be identified with known ones, but their
interpretation seems new. Using that the q-Hermite polynomials are 
closely related to the Demazure level-one characters in the twisted 
case (Sanderson, Ion), we outline a connection of our formulas to 
the level-one integrable Kac-Moody modules and coset theory.
Several instances of the level-rank duality are provided.
}

The theory of Fourier 
transform of the nilpotent double affine Hecke algebras, 
{\em Nil-DAHA}, is applied in the paper to
{\em Rogers- Ramanujan type identities\,}. Our $q$\~series  
are always modular functions (of weight zero) for the congruence 
subgroups of $SL(2,\Z)$, which has important implications. 
There are connections with the algebraic and combinatorial
theory of such identities, dilogarithm and {\em Nahm's conjecture\,},
$Y$\~systems, {\em Demazure characters\,}, the level-rank
duality and {\em coset algebras\,}. The latter are associated with 
tensor products of integrable level-one Kac-Moody representations and 
are touched only a little in this paper. We mainly study
the algebraic and arithmetic aspects. 

We arrive at an ample family of the formulas associated with 
arbitrary (reduced twisted) irreducible affine root systems 
depending on the choices of initial level-one theta functions
(numbered by minuscule weights). Some of our formulas for $p=2,3$ 
can be identified with known Rogers-Ramanujan identities, but 
there are new aspects even in these cases.
The flexibility with picking the theta functions results in the 
restricted summations in our identities (for instance, the sums
can be even or odd for $A_1$). For $A_n$, there are $\binom{p+n}{n}$ 
such choices at the level $p$ (the number of theta functions in the 
product). 

One of the main messages of the paper is that
{\em nil-DAHA can be used to calculate the key string functions for
Kac-Moody algebras of type $A_n$} \cite{KP} through the expected 
level-rank duality.

\smallskip

The following topics seem inevitable to put the results of this
paper into perspective, but it will be not a systematic review 
and only basic references will be given. We try to stick 
here and in the paper mostly to (relatively) known
examples; the level-two formulas are the major particular cases 
we consider. Quite a few
topics discussed in the introduction are not touched upon
in the main body of this work, so the following is focused on
motivation and links to other theories. 

\subsection{Using q-Hermite polynomials}
The reproducing kernel of the Fourier transform of the Nil-DAHA,
the {\em global $q$\~Whittaker function} from \cite{ChW},
is actually the key, though it is not introduced
and needed in this particular paper. Its explicit expression 
is equivalent to knowing  
$$
CT\,(P_a(X) P_b(X) \theta(X)\,\mu(X)\,), \where 
\mu\equal\prod_{\tal=[\al,j]>0}(1-q^j X_{\al}),
$$
for all pairs of {\em $q$\~Hermite polynomials} $P_a$ (their indices 
$\,a,b\,$ are anti-dominant weights in this paper) and level-one 
theta functions $\theta$. Here the product is over all positive 
affine roots $\tal$,\, $X_{[\al,j]}\equal X_{\al}q^j$,\,
$X_{a+b}=X_aX_b$ for all weights $\,a,b\,$.
By $CT(\cdot)$, we mean the constant term
of a Laurent series in terms of $X_a$.
\smallskip 

This is what we really needed. These formulas, 
used inductively, provide expansions of arbitrary products 
of level-one theta functions in terms of the $q$\~Hermite 
polynomials, which is the essence of this paper.

Then we apply $CT(\,\cdot\,P_c\,\mu)$ to the expansions
of products of theta functions for minuscule weights $c$,
adjusting them to make the constant 
term nonzero. It results in 
the Rogers-Ramanujan type 
summations, which can be then compared with
the formulas for our products in terms of the standard 
theta functions of level $p$, say from Proposition 3.14
from \cite{KP}, or with those obtained by 
direct expansions of such products, or with various formulas
in CFT, Conformal Field Theory.

The $q$\~Hermite polynomials 
completely disappear from the resulting identities. However their 
norms, products of $q$\~factorials in the denominators, and 
special quadratic forms in the powers of $q$ in the 
numerators clearly hint on their presence in the theory.  
\smallskip

Interestingly, Rogers' proof of the celebrated Rogers-Ramanujan 
identities was also based on the $q$\~Hermite polynomials. 
The generating function for these polynomials and the formula 
connecting those for $q$ and $q^{-1}$ were the main ingredients.
See, e.g. \cite{GIS} for the modern reproduction of his
method and some its generalizations. In contrast to
his approach, our proof is linked to the global
$q$\~Whittaker function, a kind of generating function for
the $q$\~Hermite polynomials of $q$\~quadratic type, instead 
of the usual generating function (cf. \cite{Sus}).

Nil-DAHA can be generally used to manage the standard
generating functions (though the formulas are more involved than 
those for $A_1$), as well as the connection 
$q\leftrightarrow q^{-1}$. Thus the 
original Rogers method can be potentially extended to general
root systems, but we do not discuss it in this work. We note that 
quite a few multivariable Rogers-Ramanujan type identities in 
the literature are actually of rank one (for 
$A_1$) but for higher levels. Here one can (and is supposed to)
use the generating function and other special features of the 
classical rank-one $q$\~Hermite polynomials at full potential.

Using $q$\~Hermite polynomials defined for arbitrary root systems 
seems new in the theory of Rogers-Ramanujan identities.
However, let us mention formula (1.2) from \cite{Ter},
where one can see a reduction of our expansion formulas
in terms of the $q$\~Hermite polynomials to a single variable
in the case of $D_n$. Let us also mention \cite{BCKL};
it is not directly connected with our approach, but
the limits of the Macdonald polynomials appear there.
\smallskip

We note that our approach generally results in $q$\~series where
the quadratic forms are multiplied (sometimes divided) 
by $2$ versus the ``main stream" of Rogers-Ramanujan type identities,
which does not mean of course that they are brand new. For instance,
our identities can be identified with known ones for the classical 
root systems and level two; what is new is their interpretation in 
terms of the $q$\~Hermite
polynomials and (sometimes) their modular invariance. We note
that this is not 
always clear in what sense known families of the Rogers-Ramanujan 
type identities are associated with root systems, even if they
look as such. For instance, the identities from \cite{An,War} 
seem associated with $B,C$, but this is not a formal link. 
\smallskip

We mention that our procedure is not that smooth
in the $q,t$\~theory based on the Macdonald polynomials 
instead of the $q$\~Hermite ones. We still can obtain interesting 
identities, but the values of the Macdonald polynomials will be 
present in these identities, and in a significant way. This 
is unwanted, since not much is known about the meaning of the
values of the Macdonald polynomials beyond the evaluation formula 
and some explicit formulas in lower ranks. Also, the $q$\~positivity 
of (all) our formulas and those for the coefficients 
of $q$\~Hermite polynomials has no known counterpart in 
the $q,t$\~ theory.
\smallskip

\subsection{Dual Demazure characters}
The $q$\~positivity mentioned above is directly related to the 
geometric meaning of the $q$\~Hermite polynomials and
the interpretation of our construction in terms of 
coset theory. The second direction is touched upon only
a little in this paper; it will be hopefully developed in other
works.  We certainly cannot interpret at the moment all 
Rogers-Ramanujan type identities we can obtain via Nil-DAHA in 
terms of Kac-Moody representation theory,
more specifically, via various coset algebras. However, the key 
step, which is the expansion of the level-one theta function 
in terms of the $q$\~Hermite polynomials, can be explained 
(though this is not a theorem at the moment). It is basically as follows.
\smallskip

Let $M$ be an irreducible level-one integrable module of 
the Kac-Moody algebra $\hat{\mathfrak{g}}$ in the simply-laced case,
$v_a$ its highest weight vector of weight $a$ with respect
to the Borel subalgebra $\hat{\mathfrak{b}}_+$.
By the {\em dual Demazure filtration} of $M$, we mean
$\{\f_b\equal U(\hat{\mathfrak{b}}_-)\,v_b\}$ for dominant $b$ 
(for $\hat{\mathfrak{b}}_+$) from
the orbit of $a$ under the action of the (extended) 
affine Weyl group, where $v_b$ is the corresponding 
{\,\em extremal vector\,} $v_b\in M$. 
Consider the corresponding adjoint (graded) module
$$
M^{ad}\ =\oplus_{b}\,\f_b^{ad},\, \where \f_b^{ad}=
\f_b/(+_{c>b\,}\f_c) 
\hbox{\ \,for dominant\ \,} b,c
$$  
in terms of the standard ordering of dominant weights
(which are partitions for $A_n$).  Then 
the modules $\f_b^{ad}$ can be identified with the
so-called {\em global Weyl modules\,}; see \cite{FeL} and
\cite{FoL}. 

The claim is that the character of $\,\f_b^{ad}\,$
is the character of the corresponding 
{\em local Weyl module\,} upon its multiplication by $q^{b^2/2}$ 
and division by the product $\prod_{i=1}^n\prod_{j=1}^{m_i} (1-q^j)$, 
where $b=\sum m_i\om_i$ for fundamental $\{\om_i\}$.
Furthermore, the character of the local Weyl module here is
the {\em dual Demazure character} defined as  
$D^\star_b(X,q)\,=\,D_b(X^{-1},q^{-1})$ for
the standard Demazure character $D_b$. Here the
substitution $X_a=e^{-a}$ establishes the connection
with the standard notation,
the character is the trace of $q^{L_0}$
for the energy operator  $L_0$ from the Virasoro algebra. 

The connection of the dual Demazure filtration with the
Weyl modules is known in the simply-laced case 
(see, e.g. \cite{FoL}), but we can not give an exact reference
concerning the formula for their characters  
in terms of the Demazure characters. Though, see formula (3.25) 
from \cite{FJKMT} in the case of $A_1$.  Such relation 
seems not fully established at the moment. We note that using 
the DAHA-based identities from this paper can be used for the 
justification of this relation and similar facts. Indeed, generally 
we know that the formulas in terms of the Demazure operators give
characters of modules {\em no smaller} than the actual ones; then 
we can use the identities obtained in this paper.  

The definition above results in the
equality  $D^\star_b(X,q)=P_{-b}(X,q)$
for the $q$\~Hermite polynomials $P_{-b}$ due to
\cite{San} ($G\!L$) and \cite{Ion1}
(arbitrary reduced root systems).
They established that the level-one Demazure characters in the 
twisted case are 
$P_{-b}(X^{-1},q^{-1})$ using the DAHA-intertwiners. The
substitution $q\mapsto q^{-1}$ is important here. Changing
$X$ to $X^{-1}$ is not too significant for dominant weights $\,b\,$,
but this connection holds for any nonsymmetric $q$\~Hermite 
polynomials.  
Let us mention here that the modified Hall-Littlewood 
polynomials, known to be connected to Demazure characters in
some examples, are  closely related to the $q$\~Hermite polynomials.
\smallskip

\subsection{Toward coset models}\label{Tcosets}
The representation theory interpretation
of our identities will be subject of our further research. 
However, it is important to explain in this
paper why we think that the family of Rogers-Ramanujan
type identities we obtain is essentially in one-to-one 
correspondence with the coset 
decomposition of tensor products of level-one representations, 
certainly one of the key problems in coset theory. 

We refer
to \cite{Kac,Kum} for the necessary definitions; also, see
paper \cite{FJMT}, especially
formulas (1.4)-(1.6) there, devoted to the
matters closely related to what we discuss below.

Let $M_1,M_2,\ldots, M_p$ be a collection of irreducible integrable 
representations of $\hat{\mathfrak{g}}$ of levels
$l_1,\ldots,l_p$\,, $L=L_{\{\la,l\}}$ an irreducible integrable 
representation of level $\,l=l_1+\cdots l_p$ with
the highest weight $\{\la,l\}$ for a (nonaffine) dominant 
weight $\la$.  We consider the highest weight modules with respect 
to $\hat{\mathfrak{b}}_+$. Let
\begin{align*}
&\nu(L\,;\, M_1,\cdots, M_p)\ \equal\  
\hbox{Hom}_{\hat{\mathfrak{g}}}\,
(L\,, M_1\otimes M_2\otimes\ldots\otimes M_p).
\end{align*}

This space can be expected to
be an irreducible module of the {\em coset vertex operator
algebra} defined essentially as the commutant (centralizer) 
of $U(\hat{\mathfrak{g}})$ diagonally embedded into 
$U(\hat{\mathfrak{g}}\times\ldots\times \hat{\mathfrak{g}})$
\,($p$ times).
Importantly, the Virasoro algebra belongs to
the coset algebra; using the energy operator $L_0$ we set
\begin{align*}
&\chi_q(M_1,\ldots, M_p) \and \chi_q(L\,;\,M_1,\ldots, M_p)\ =\ 
\hbox{trace\,} (\,q^{L_0}\,) \\
&\hbox{for\ } L_0 \hbox{\ acting in\ }
M_1\otimes M_2\otimes\ldots\otimes M_p,\,
\nu(L\,;\, M_1,\ldots, M_p).
\end{align*}

Let $\hat{\mathfrak{b}}_-$ be the Borel subalgebra
opposite to $\hat{\mathfrak{b}}_+$,  
$\mathfrak{h}$ the Cartan subalgebra and
$\C_{-\la}$ the one dimensional 
$\hat{\mathfrak{b}}_-$\~module of weight
$-\la$ (for the weight $\la$ above).  A standard fact in Kac-Moody 
theory is that
\begin{align*}
&\nu(L\,;\, M_1,\ldots, M_p)\ =\  
H_\star (\hat{\mathfrak{b}}_-,\mathfrak{h}\,;\, 
M_1\otimes M_2\otimes\ldots\otimes M_p\otimes \C_{-\la}),\\
&=\ H_0 (\hat{\mathfrak{b}}_-,\mathfrak{h}\,;\, 
M_1\otimes M_2\otimes\ldots\otimes M_p\otimes \C_{-\la}).
\end{align*}
Here $H_\star (\hat{\mathfrak{b}}_-,\mathfrak{h}\,;\,\cdot\,)$
is relative homology.
The higher homology vanishes due to the
integrability of the modules we consider;
see, e.g. \cite{Kum}.  

Thus we can identify 
$\chi_q(L\,;\,M_1,\ldots, M_p)$
with the Euler characteristic of the complex 
$$
\bigl[\,\wedge^\star(\hat{\mathfrak{h}}_-)\otimes 
M_1\otimes\ldots\otimes M_p\otimes \C_{-\la}\,\bigr]^\mathfrak{h},
$$
which, in turn, is 
$
\,CT\,(\,\chi_q(M_1)\cdot\ldots\cdot\chi(M_p)\,
\chi_q(\C_{-\la})\,\mu\,).$ 
When $c=0$, i.e. for the vacuum representation
$L$ of level $l$, it will be the constant term of 
$
\chi_q(M_1)\cdot\ldots\cdot\chi(M_p)\,\mu\,;
$ 
the character $\chi_q(\C_{-\la})$ is $q^{-\la}$, which
is $X_{-\la}$ in our notation. 

Alternatively,  one can express 
$M_1\otimes\ldots\otimes M_p$ as a sum of
irreducible integrable $\hat{\mathfrak g}$\~modules and then 
use the Weyl-Kac character formula for each of them; 
$\mu$ is essentially the denominator of this formula.



We see that our identities can be used to determine
the characters of the coset algebra acting in   
the spaces $\nu(L\,;\, M_1,\ldots, M_p)$
for level-one $M_1,\ldots, M_p$. The other way around,
one can use these characters (when they are known)
to obtain interesting expressions for the Rogers-Ramanujan
type series from this paper.
\smallskip
 
\subsection{Around Nahm's conjecture}
For arbitrary root systems and levels, we arrive at
$q$\~series in the form
$$
F_{A,B,C}(q)\ =\ 
\sum_{n\in \Z_+^r}\frac{q^{n^{T}An/2+Bn+C}}{(q_1)_{n_1}
\cdots(q_r)_{n_r}},\ \, (q)_m=\prod_{i=1}^m (1-q^i).
$$
For $p=2$, $A$ is the inverse Cartan matrix ($r$ is the rank),
$q_i=q^{\nu_i}$, where  $\nu_i=1$ for short simple roots $\al_i$
and $\nu_i=2,3$ for long simple $\al_i$  correspondingly
for $B_n$-$C_n$-$F_4$ and $G_2$. In the simply-laced case,
it is exactly the class of series from
the so-called {\em Nahm's conjecture} \cite{Na}, generally, 
for symmetric real positive definite matrices $A$.  

See \cite{KM,KN,NRT,KKMM} concerning the physics origins of this 
conjecture, they are not far from our approach (related to 
the Verlinde algebras). 
We will not discuss here  the key role of the 
thermodynamic Bethe ansatz (TBA) due to Al. Zamolodchikov
and others,
the $Y$\~systems and the cluster algebras; 
see, e.g. \cite{IKNS,Nak,Nak1} for recent developments.
Also, applications of dilogarithms to volumes of
$3$\~manifolds will be completely omitted. 

We note that the method from \cite{NRT,Na} does not 
rely on TBA and is actually similar to that from \cite{RS}
(direct calculating the saddle point). The formula for the 
sum of $L(Q_i)$  (the next section) was obtain in \cite{RS}
for the inverse Cartan matrix of type $A_n$  using 
the asymptotic formula for the partition function $p(k)$;
another proof is in \cite{KR}.

The simplest cases of our formulas for  $A_1$ and levels $p=2,3$
for the unrestricted theta function (no parity constraints in
the summation) can be found in Tables 1,2 from \cite{Za}; 
see also Theorems 3.3 and Theorem 3.4 from \cite{VZ}. For
$p=3$, we found some new developments.

When the root systems is $A_2$ and $p=2$, the matrix $A$ is
{\tiny $\left(
     \begin{array}{cc}4/3 & 2/3 \\2/3 & 4/3 \\ \end{array}\right).
$}
In this case, our formula is also from \cite{Za},\cite{VZ}, 
but we can obtain more identities (generally many more) by using 
various (level-one) theta functions in the products.

Our matrix $A$ is associated with the standard bilinear form
in tensor product of the weight lattice $P_R$
(any reduced root system $R$, the twisted case) and the root 
lattice $Q$ of type $A_{p-1}$; $p$ is the level. Such class of 
$A$ attracts a lot of attention in related physics. 
The summation can be restricted in our formulas by picking any 
decomposition of $p$ as a sum of $|P_R/Q_R|$ nonnegative terms 
counting the numbers of theta functions in the $p$\~product associated 
with the corresponding minuscule weights. For instance,
the number of such choices for $R=A_{n}$ equals the number
of decompositions $p=a_1+\cdots+a_{n+1}$, where $a_i\ge 0$
and the order matters; thus it equals $\binom{n+p}{n}$.
\smallskip

All our series are 
modular functions. It allows verifying some
conjectural Rogers-Ramanujan type identities and finding new ones.  
For instance, 
among other applications, our approach provides a
justification of two formulas with question marks 
in Table 1 from \cite{VZ} (for us, this is the case of $A_2$, $p=2$).
Also we found a split of the formulas from 
Theorem 3.4 there for $A=$
{\tiny $\left(
     \begin{array}{cc}1 & -1/2 \\-1/2 & 1 \\ \end{array}\right)
$},
which relates them to the classical 
Rogers-Ramanujan identities. It requires using the parity 
restrictions in the summation, our flexibility with level-one 
theta functions; this is the case of $A_1, p=3$.

We note that the Nahm's conditions hold 
for $p=2$ (all root systems), which was established in 
\cite{Lee} using $Y$\~systems;
we thank Tomoki Nakanishi for the reference. Our approach, 
combined with Corollary 3.2 from \cite{VZ} (slightly 
modified to include the $BCFG$\~types), gives the rationality 
of $\sum_i\,L(Q_i)$ (see below) for any levels $p$, 
but only for the (unique) distinguished solutions in the 
range $\{0<Q_i<1\}$.   
\medskip

\comment{
Let us mention here that the affine 
Toda field theory is conjecturally connected with the $S$\~matrices 
associated with $\{Q_i\}$; this is discussed in \cite{KM} but without
much detail. The $q$\~Hermite polynomials are closely related
to the {\em multi-dimensional\,} 
$q$\~Toda operators, which in their turn are connected
with the affine (periodic) Toda chains.

It is worth mentioning here that the {\em diagonal\,} factorizable 
$S$\~matrices, the foundation of the TBA, are closely related to 
the $S$\~matrices describing factorizable particles on a line and 
satisfying the (nondiagonal) Yang-Baxter equations.  The latter 
direction, started by Zamolodchikovs, was extended by the first 
author to arbitrary root systems (including the affine ones). 
The classical $BC$ root 
systems were associated with the particles on a half-line (with the
reflection at the endpoint) and on a segment. It is directly connected
with DAHA via the quantum Knizhnik-Zamolodchikov equations.
It is interesting to reexamine these links from the viewpoint
of the present paper.
}

To conclude this general discussion (more detail will be provided
in the next section), we must note that not many formulas listed
in \cite{Za,VZ} and in the vast literature on the Rogers-Ramanujan
identities can be obtained by our construction. Recall that only 
modular invariant $q$\~series result from our approach, but certainly
not all such. Hopefully, 
using the root system $C^\vee C_n$ and $l\in \Z/2$ 
(see the end of the paper) will significantly increase the scope 
of our approach. Also, higher levels $p$ can be generally used 
followed by various reductions, but this possibility was not
systematically explored at the moment. 
\medskip

\subsection{Dilogarithm identities}
Continuing the previous section, let as briefly discuss
using the Rogers dilogarithm 
$$
L(z)\ \equal\ Li_2(z)+\frac{1}{2}\log x \log(1-x),
$$
the key in Nahm's conjecture. The modular invariance of the
$q$\~series $F_{A,B,C}(q)$ 
above implies the following (see \cite{VZ},
\cite{Za} and \cite{Na}).

For a $N\times N$\~matrix $A=(a_{ij})$ 
the system of equations
\begin{align*}
1-Q_i=\prod_{j=1}^N Q_j^{a_{ij}},\ i=1,\ldots, N
\end{align*}
has a unique solution in the range $0<Q_i<1$,
assuming that $A$ is real symmetric and positive 
definite. It is referred to as well-known in 
physics papers; this system
is part of the Thermodynamic Bethe Ansatz (TBA).
See Lemma 2.1 from \cite{VZ} for a direct 
justification.

The main claim is that 
$$
L_{A}\ =\ \frac{6}{\pi^2}\,\sum_{i=1}^N L(Q_i)\ \in \ \Q,
$$
assuming the modular invariance of $F_{A,B,C}(q)$,
which can be interpreted as finding a torsion element in 
the corresponding Bloch group. See \cite{Za} and \cite{VZ}
for the justification.
Nahm's conjecture states that any (complex) solutions  
$\{Q_i\}$ result in torsion elements in the Bloch group
in the same way.
Moreover, the latter property is equivalent to 
the modular invariance of $F_{A,B,C}(q)$ 
for suitable $B,C$. It appeared generally not the case
\cite{VZ}, but this conjecture certainly clarifies the 
role of dilogarithms here.

It is of obvious interest to analyze $L_A$ for our 
$q$\~series and related ones. 
See \cite{KM} and \cite{KN,NRT,KKMM,Ter} 
for the $A$-$D$-$E$ cases for $p=2$ and papers 
\cite{IKN,IKNS,Nak,Nak1} for recent developments.
The $A$\~case ($p=2$) was calculated in \cite{RS} and \cite{KR}
(see also \cite{FrS}).

We note that the tadpole case $T_n$ corresponds to $C_n$,
where $A=(a_{ij})$ for 
$a_{i,j}=2\,\hbox{Min}(i,j)\,$, which is proportional 
to the inverse of the ``Cartan matrix" of $T_n$. The
corresponding $Q$\~system is {\em exactly} that for $A_{2n}$
upon symmetry $Q_{i}=Q_{2n-i+1}\, (1\le i\le n)$.
\smallskip

Let us begin with $A_3,A_4$ and $D_4$ in the case of
level $2$. Using the uniqueness of $\{Q_i\}\subset (0,1)$,
we can impose the symmetries resulting from those
in the corresponding inner products.
In these cases, the $Q$\~systems, their solutions 
the corresponding $L=L_A$ are:
\begin{align*}
A_3:\ &1-Q_1=Q_1^{\frac{3}{2}}Q_2Q_3^{\frac{1}{2}}, 
1-Q_2=Q_1Q_2^2 Q_3,
1-Q_3=Q_1^{\frac{1}{2}}Q_2 Q_3^{\frac{3}{2}},\\
 &\hbox{setting\ }Q_1=Q_3,\ \ Q_1=2/3=Q_3,\ Q_2=3/4
\and  L=2;\\
A_4:\ &1-Q_1=Q_1^{8/5}Q_2^{6/5}Q_3^{4/5}Q_4^{2/5},\  
1-Q_2=Q_1^{6/5} Q_2^{12/5} Q_3^{8/5}Q_4^{4/5},\\
 &1-Q_3=Q_1^{4/5}Q_2^{8/5}Q_3^{12/5}Q_4^{6/5},\ 
1-Q_4=Q_1^{2/5}Q_2^{4/5}Q_3^{6/5}Q_4^{8/5},\\ 
 &\hbox{setting\ }
Q_1=Q_4,\, Q_2=Q_3, \ \ \ Q_1=1-Q_2^{-2}+Q_2^{-1} \and \\ 
&Q_2\in (0,1) \hbox {\ \,is a unique solution of\ }\, 
t^3+2t-t-1=0\,: \\
 &Q_2=2\cos(\frac{\pi}{7})\!-\!1 \hbox{ and }
L_{A_4}=\frac{20}{7}\ \hbox{(cf. Watson's identities)};\\
D_4:\ &1-Q_1=Q_1^{2}Q_2^{2}Q_3Q_4,\  
1-Q_2=Q_1^{2} Q_2^{4} Q_3^{2}Q_4^{2},\\
 &1-Q_3=Q_1 Q_2^{2}Q_3^{2}Q_4,\ 
1-Q_4=Q_1Q_2^{2}Q_3Q_4^{2},\\ 
 &\hbox{for\ }
Q_1\!=\!Q_3\!=\!Q_4,\  
Q_1=\frac{3}{4}\,,\,Q_2=\frac{8}{9}\,,\ \,L=3. 
\end{align*}

More generally, our $L$\~sums are exactly those found in
\cite{KM} times $h/2$ for the Coxeter number
$h=(\rho,\vth)+1$. Note that our $Q$\~systems must be transformed
to match \cite{KM} following formulas (63,71) there; see formulas
(73,77,79,81,83) and (A1) from
\cite{Ter} for $T_n$ (directly related to $A_{2n}$).
Namely,
\begin{align*}
&L=\frac{n(n+1)}{(n+3)} \for A_n, \ \ 
Q_n\,=\,2\cos(\frac{\pi}{n+3})-1 \when  n=2m,\\ 
&Q_n=\frac{1+\cos(\frac{\pi}{m+2})}{2} \for n=2m+1\,
(\hbox{see formulas (1,2) in \cite{RS}})\,; \\ 
&L=n-1 \for D_n(n>3), \where \\
&Q_i=\frac{(i+1)^2-1}{(i+1)^2}\for
i<n-1 \and Q_{n-1}=\frac{n-1}{n}=Q_n\,;\\
&L=\frac{n(2n+1)}{(2n+3)} \for T_n\,;\ \ 
L=\,\frac{36}{7},\,\frac{63}{10},\,\frac{15}{2} \for E_{6,7,8}.
\end{align*}

The $Q_i$ are well-known in the $A_n$\~case; see \cite{KM}
or \cite{RS,KR,KN,FrS}.
In physics and RT literature,  $2L/h=c_{\hbox{\tiny eff}}\,$ 
is called the effective central charge
(the finite-size scaling coefficient). 
It is the difference $c-12\,d_0\,$ for the corresponding CFT. 
As such, this is of importance to us, since it provides information
on the structure of the right-hand side of the Rogers-Ramanujan
type identities we obtain and because we can determine from it
which coset models can be expected. See \cite{KM}
and the end of this paper for some discussion. 

For an arbitrary root system $R\subset \R^n$, we normalize the
form by the condition 
$(\al_{\hbox{\tiny short}},\al_{\hbox{\tiny short}})=2$. 
Then our $A$ is $2$ times the restriction of $(\cdot,\cdot)$ 
to the weight lattice. We need to modify the $Q$\~system
and $L$\,:
\begin{align*}
(1-Q_i)^{\nu_i}
=\prod_{j=1}^n Q_j^{a_{ij}}\, (1\le i\le n),\ \,
L_R = \frac{6}{\pi^2}\,\sum_{i=1}^n \nu_i
L(Q_i),\  \nu_i=\frac{(\al_i,\al_i)}{2}.
\end{align*}  
Then our analysis based on \cite{VZ} results in the following:
\begin{align*}
L_{B_n}=\frac{n(2n-1)}{n+1},\ \, L_{C_n}=n,\ \,
L_{F_4}=\frac{36}{7},\ \, L_{G_2}=3.
\end{align*}
Note that  $L_{B_n}=L_{A_{2n-1}},\, L_{F_4}=L_{E_6}$.
For $C_n$, we have two kinds of 
dilogarithm identities, this one and that for $T_n$ 
(without using $\nu_i$ on the left-hand side).
As a matter of fact, the Rogers-Ramanujan series
of type $R\,$ in our approach corresponds to the $Q$\~system 
for $R^\vee$, namely, $C_n$ corresponds to $B_n$ and
$B_n$ corresponds to the $Q$\~system of type $T_n$. 
We will not go into detail. 

We are grateful to Tomoki Nakanishi for
identifying the formulas for $L$ in the $BCFG$ cases with the 
known instances of the (constant) $Y$\~systems via the so-called
folding construction; see Section 9 from \cite{IKNS}
and \cite{Nak,Nak2}. See also \cite{Nak1} and 
\cite{IKN}, Theorem 2.10.

\subsection{Obtaining ``Rogers-Ramanujan"}
Let us demonstrate what our approach gives for the classical
Rogers-Ramanujan series, which occur at the level $3$ for $A_1$ 
in our construction. Recall that the classical interpretation
of these identities from \cite{LW} was also associated with
$A_1$ at the level three, though our tools are different
and we can connect our formulas with the classical Rogers-Ramanujan 
identities only upon certain transformations (it actually goes via the
coset construction and is less direct
than in \cite{LW}). 

We begin with the product of three (level-one) theta functions 
for $A_1$, correspondingly even, even and odd, where
$$
\th_k(X)=\sum_{j=-\infty}^{\infty} q^{(2j+k)^2/4}X^j \for 
k=0\, (\hbox{even}),\ 1\, (\hbox{odd}).
$$

The $\mu$\~function is the classical Jacobi theta function\,:
\begin{align*}
&\mu =  
\prod_{j=0}^\infty (1-X^2 q^{j})
(1-X^{-2}q^{j+1}),\ \, 
CT(\mu) = 
\prod_{j=1}^{\infty} \frac{1}
{1-q^j}\,.
\end{align*}
The bilinear form $CT(f\,g\,\mu)$ in the space 
of symmetric Laurent polynomials in terms 
of $X^{\pm 1}$ makes the $q$\~Hermite polynomials pairwise orthogonal. 

The example of $p=3$ is related to the Rogers-Ramanujan summations
as follows. For $k=0,1$,
\begin{align*}
&\frac{CT(\, \th_1\th_0\th_k (X+X^{-1})^{1-k}\mu\,)}
{q^{\frac{1+k}{4}}\prod_{j=1}^\infty(1-q^j)^2}\ =\ 
\sum_{n,m\ge 0}\frac{q^{2(n^2-nm+m^2)+nk+m}}
{\prod_{j=1}^{2n+1}(1-q^j)\prod_{j=1}^{2m+1}(1-q^j)}\\
&=\!\!\!\sum_{n,m\ge 0}\frac{q^{2(n^2-nm+m^2)+2nk-m}}
{\prod_{j=1}^{2n+k}(1-q^j)\prod_{j=1}^{2m}(1-q^j)}\!=\!
\sum_{n=0}^\infty \frac{q^{2n^2+2nk}}{\prod_{j=1}^n(1-q^{2j})}
\prod_{j=1}^{\infty}(1+q^j)^2.
\end{align*}

It is a combination of (\ref{Rog-Ram1},\ref{Rog-Ram2})
with formulas (\ref{r-r51},\ref{r-r52}),
where we consider the sequences $1\!\0\!0$,
$0\!1\!0$ for $k=0$
and $1\!0\!1$, $1\!1\!0$ for $k=1$ respectively. Such permutations
obviously do not influence the left-hand side, but result
in different double summations.

To be more exact, the first two equalities here and their
generalizations to arbitrary ranks and levels are the main 
output of this paper. The reduction to a single
summation is special for this particular example.
See Section \ref{sec:casep=3} for 
formulas for the remaining combinations of the
indices of theta functions and further details.

We arrive at the Rogers-Ramanujan series
with $q^2$ instead of $q$\,: 
\begin{align*}
&\frac{CT(\, \th_1\th_0\th_k (X+X^{-1})^{1-k}\mu\,)}
{q^{\frac{1+k}{4}}\prod_{j=1}^\infty(1-q^{2j})^2}\ =\ 
\sum_{n=0}^\infty \frac{q^{2n^2+2nk}}{\prod_{j=1}^n(1-q^{2j})}.
\end{align*}

\comment{
\begin{align*}
&\frac{CT(\, 
\th_1\th_0\th_k (X+X^{-1})^{1-k}\mu\,)}
{q^{\frac{1+k}{4}}\prod_{j=1}^\infty(1-q^j)^2}\ =\\
&\sum_{n,m\ge 0}\frac{q^{2(n^2-nm+m^2)+nk+m}}
{\prod_{j=1}^{2n+1}(1-q^j)\prod_{j=1}^{2m+1}(1-q^j)}
=
\sum_{n=0}^\infty \frac{q^{2n^2+2nk}}{\prod_{j=1}^n(1-q^{2j})}
\prod_{j=1}^{\infty}(1+q^j)^2.
\end{align*}
\begin{align*}
&q^{-1/4}\frac{CT(\, \th_1\th_0^2 (X+X^{-1})\mu\,)}
{\prod_{j=1}^\infty(1-q^j)^2}\ =\\
&\sum_{n,m\ge 0}\frac{q^{2(n^2-nm+m^2)-m}}
{\prod_{j=1}^{2n}(1-q^j)\prod_{j=1}^{2m}(1-q^j)}
=
\sum_{n=0}^\infty \frac{q^{2n^2}}{\prod_{j=1}^n(1-q^{2j})}
\prod_{j=1}^{\infty}(1+q^j)^2.
\end{align*}
\begin{align*}
&q^{-1/4}CT(\, \th_1\th_0^2 (X+X^{-1})\mu\,)
=
\sum_{n=0}^\infty \frac{q^{2n^2}}{\prod_{j=1}^n(1-q^{2j})}
\prod_{j=1}^{\infty}(1-q^{2j})^2.
\end{align*}
}

At the end of this paper an outline of the interpretation
of this formula for $k=1$ in terms of coset theory
is provided and we establish a relation of these formulas
to the string functions for $A_2$ of level two.
\medskip

{\bf Acknowledgements.} We would like to thank the
Mathematics Department of the Kyoto University and RIMS 
for the invitations and hospitality; our special thanks
to Tetsuji Miwa. We thank Sergei Loktev for a useful discussion
and George Andrews for important suggestions. 
Our special thanks to Ole Warnaar for attracting our attention
to \cite{War}, useful comments on the paper and valuable
help with the level-two $B,C$\~series in Section \ref{sec:lev2BC}.
We are very grateful to Tomoki Nakanishi for identifying the
dilogarithm identities we arrive at with proper instances of
the constant $Y$\~systems; see \cite{Nak2}.
\medskip

\setcounter{equation}{0}
\section{Theta-products via Nil-DAHA}
Let $R=\{\al\}   \subset \R^n$ be a root system of type
$A,B,...,F,G$
with respect to a euclidean form $(z,z')$ on $\R^n
\ni z,z'$,\ 
$W$ the {\em Weyl group} 
generated by the reflections $s_\al$,\ 
$R_{+}$ the set of positive  roots
corresponding to fixed simple 
roots $\al_1,...,\al_n,$\ 
$\Ga$ the Dynkin diagram
with $\{\al_i, 1 \le i \le n\}$ as the vertices,\ 
$\rho=\frac{1}{2}\sum_{\al\in R_+} \al$,\ 
$R^\vee=\{\al^\vee =2\al/(\al,\al)\}.$

The root lattice and the weight lattice are:
\begin{align}
& Q=\oplus^n_{i=1}\Z \al_i \subset P=\oplus^n_{i=1}\Z \om_i,
\notag
\end{align}
where $\{\om_i\}$ are fundamental weights:
$ (\om_i,\al_j^\vee)=\de_{ij}$ for the
simple coroots $\al_i^\vee.$
Replacing $\Z$ by $\Z_{\pm}=\{k\in\Z, \pm k\ge 0\}$ we obtain
$Q_\pm, P_\pm.$
Here and further see  \cite{Bo}.

The form will be normalized
by the condition  $(\al,\al)=2$ for the
{\em short} roots in this paper. 
Thus
$\nu_\al\equal (\al,\al)/2$ can be either $1,$ $\{1,2\},$ or
$\{1,3\}.$ 
\noindent
This normalization leads to the inclusions
$Q\subset Q^\vee,  P\subset P^\vee,$ where $P^\vee$ is
defined to be generated by
the fundamental coweights $\{\om_i^\vee\}$ dual to $\{\al_i\}$. 

We note that $Q^\vee=P$ for $C_n (n\ge 2)$, $P\subset Q^\vee$
for $B_{2n}$ and $P\cap Q^\vee=Q$ 
for $B_{2n+1}$; the index $[Q^\vee:P]$
is $2^{n-2}$ for any $B_n$ (in the
sense of lattices).

\subsection{Affine Weyl groups}
The vectors $\ \tal=[\al,\nu_\al j] \in
\R^n\times \R \subset \R^{n+1}$
for $\al \in R, j \in \Z $ form the
{\em affine root system}
$\tR \supset R$; this is the so-called twisted case.
\smallskip

The vectors $z\in \R^n$ are identified with $ [z,0]$.
We add $\al_0 \equal [-\vth,1]$ to the simple
roots for the {\em maximal short root} $\vth\in R_+$.
It is also the {\em maximal positive
coroot} because of the choice of normalization. 
The Coxeter number is then $h=(\rho,\vth)+1$.
The corresponding set
$\tR_+$ of positive roots equals 
$R_+\cup \{[\al,\nu_\al j],\ \al\in R, \ j > 0\}$.

We complete the Dynkin diagram $\Ga$ of $R$
by $\al_0$ (by $-\vth$, to be more
exact); this is called the {\em affine Dynkin diagram}
$\tGa$. One can obtain it from the
completed Dynkin diagram from \cite{Bo} for 
the {\em dual system}
$R^\vee$ by reversing all arrows. 

The set of
the indices of the images of $\al_0$ by all
the automorphisms of $\tGa$ will be denoted by $O$;\,
$O=\{0\} \for E_8,F_4,G_2$. Let $O'\equal \{r\in O, r\neq 0\}$.
The elements $\om_r$ for $r\in O'$ are the so-called minuscule
weights: $(\om_r,\al^\vee)\le 1$ for
$\al \in R_+$ (here  $(\om_r,\vth)\le 1$ is sufficient).
\smallskip

{\bf Extended Weyl groups.}
Given $\tal=[\al,\nu_\al j]\in \tR,  \ b \in P$, the
corresponding reflection in $\R^{n+1}$ is defined by the formula
\begin{align}
&s_{\tal}(\tz)\ =\  \tz-(z,\al^\vee)\tal,\
\ b'(\tz)\ =\ [z,\ze-(z,b)],
\label{ondon}
\end{align}
where $\tz=[z,\ze] \in \R^{n+1}$.
\smallskip

The {\em affine Weyl group} $\tW$ is generated by all $s_{\tal}$
(we write $\tW = \lan s_{\tal}, \tal\in \tR_+\ran)$. One can take
the simple reflections $s_i=s_{\al_i}\ (0 \le i \le n)$
as its
generators and introduce the corresponding notion of the
length. This group is
the semidirect product $W\lsmash Q'$ of
its subgroups $W=$ $\lan s_\al,
\al \in R_+\ran$ and $Q'=\{a', a\in Q\}$, where
\begin{align}
& \al'=\ s_{\al}s_{[\al,\,\nu_{\al}]}=\
s_{[-\al,\,\nu_\al]}s_{\al}\for
\al\in R.
\label{ondtwo}
\end{align}

The {\em extended Weyl group} $ \hW$ generated by $W\and P'$
(instead of $Q'$) is isomorphic to $W\lsmash P'$:
\begin{align}
&(wb')([z,\ze])\ =\ [w(z),\ze-(z,b)] \for w\in W, b\in B.
\label{ondthr}
\end{align}
From now on  $b$ and $b',$ $P$ and $P'$ will be identified.

Given $b\in P_+$, let $w^b_0$ be the longest element
in the subgroup $W_0^{b}\subset W$ of the elements
preserving $b$. This subgroup is generated by simple
reflections. We set
\begin{align}
&u_{b} = w_0w^b_0  \in  W,\ \pi_{b} =
b( u_{b})^{-1}
\ \in \ \hW, \  u_i= u_{\om_i},\pi_i=\pi_{\om_i},
\label{xwo}
\end{align}
where $w_0$ is the longest element in $W,$
$1\le i\le n.$

The elements $\pi_r\equal\pi_{\om_r}, r \in O'$ and
$\pi_0=\hbox{id}$ leave $\tGa$ invariant
and form a group denoted by $\Pi$,
 which is isomorphic to $P/Q$ by the natural
projection $\{\om_r \mapsto \pi_r\}$. As to $\{ u_r\}$,
they preserve the set $\{-\vth,\al_i, i>0\}$.
The relations $\pi_r(\al_0)= \al_r= ( u_r)^{-1}(-\vth)$
distinguish the
indices $r \in O'$. Moreover,
\begin{align}
& \hW  = \Pi \lsmash \tW, \where
  \pi_rs_i\pi_r^{-1}  =  s_j \iif \pi_r(\al_i)=\al_j,\
 0\le j\le n.
\end{align}

{\bf The length.}
Setting
$\hw = \pi_r\tw \in \hW,\ \pi_r\in \Pi, \tw\in \tW,$
{\em the length} $l(\hw)$
is by definition the length of the reduced decomposition
$\tw= $ $s_{i_l}...s_{i_2} s_{i_1} $
in terms of the simple reflections
$s_i, 0\le i\le n.$ The number of  $s_{i}$
in this decomposition
such that $\nu_i=\nu$ is denoted by   $l_\nu(\hw).$

The length can be
also defined as the
cardinality $|\la(\hw)|$
of the {\em $\la$\~set} of $\hw$\,:
\begin{align}\label{lasetdef}
&\la(\hw)\equal\tR_+\cap \hw^{-1}(\tR_-)=\{\tal\in \tR_+,\
\hw(\tal)\in \tR_-\},\
\hw\in \hW.
\end{align}
Alternatively,
\begin{align}
&\la(\hw)=\cup_\nu\la_\nu(\hw),\
\la_\nu(\hw)\
\equal\ \{\tal\in \la(\hw),\nu({\tal})=\nu \}.
\label{xlambda}
\end{align}
See, e.g. \cite{Bo,Hu} and also \cite{C4,C101}.

\subsection{Nil-DAHA}
For pairwise commutative $X_1,\ldots,X_n,$ let
\begin{align}
& X_{\tb}\ =\ \prod_{i=1}^nX_i^{l_i} q^{k}
\iif \tb=[b,k],\ \hw(X_{\tb})\ =\ X_{\hw(\tb)}.
\label{Xdex}\\
&\hbox{where\ } b=\sum_{i=1}^n l_i \om_i\in P,\ j \in
\Q,\ \hw\in \hW.
\notag \end{align}
For instance, $X_0\equal X_{\al_0}=qX_\vth^{-1}$.
We will set $(\tilde{b},\tilde{c})=(b,c)$, ignoring the 
affine extensions in this pairing.

Note that $\pi_r^{-1}$ is $\pi_{r^*}$ and
$u_r^{-1}$ is $u_{r^*}$
for $r^*\in O\ ,$ where  
the reflection $^*$ is
induced by an involution of the nonaffine Dynkin diagram
$\Gamma.$ 
By  $m,$ we denote the least natural number
such that  $(P,P)=(1/m)\Z.$  Thus
$m=2 \for D_{2k},\ m=1 \for B_{2k} \and C_{k},$
otherwise $m=|\Pi|$.

\begin{definition}
The nil-DAHA \,$\HH\, $
is generated over $\Z_q\equal \Z[q^{\pm 1/m}]$ by
the elements $\{ T_i,\ 0\le i\le n\}$,
pairwise commutative $\{X_b, \ b\in P\}$ satisfying
(\ref{Xdex}),
and the group $\Pi,$ where the following relations are imposed:

(o)\ \,\,\,$ T_i(T_i+1)\ =\
0,\ 0\ \le\ i\ \le\ n$;

(i)\ \ \ $ T_iT_jT_i...\ =\ T_jT_iT_j...,\ m_{ij}$
factors on each side;

(ii)\ \   $ \pi_rT_i\pi_r^{-1}\ =\ T_j \iif
\pi_r(\al_i)=\al_j$;

(iii)\  $T_iX_b \ =\ X_b X_{\al_i}^{-1} (T_i+1)\} \iif
(b,\al^\vee_i)=1,\
0 \le i\le  n$;

(iv)\ $T_iX_b\ =\ X_b T_i$ if $(b,\al^\vee_i)=0
\for 0 \le i\le  n$;

(v)\ \ $\pi_rX_b \pi_r^{-1}\ =\ X_{\pi_r(b)}\ =\
X_{ u^{-1}_r(b)}
 q^{(\om_{r^*},b)},\  r\in O'$.
\label{double}
\end{definition}

{\bf T-elements.}
Note that one can rewrite (iii,iv) as in \cite{L}:
\begin{align}
&T_iX_b -X_{s_i(b)}T_i\ =\
\frac{X_{s_i(b)}-X_b}
{1-X_{\al_i}},\ 0 \le i\le  n.
\label{tixi}
\end{align}

Given $\tw \in \tW, r\in O,\ $ the product
\begin{align}
&T_{\pi_r\tw}\equal \pi_r\prod_{k=1}^l T_{i_k},\where
\tw=\prod_{k=1}^l s_{i_k},
l=l(\tw),
\label{Twx}
\end{align}
does not depend on the choice of the reduced decomposition
(because $T_i$ satisfy the same ``braid'' relations
as $s_i$ do).
Moreover,
\begin{align}
&T_{\hv}T_{\hw}\ =\ T_{\hv\hw}\  \hbox{ whenever}\
 l(\hv\hw)=l(\hv)+l(\hw) \for
\hv,\hw \in \hW. \label{TTx}
\end{align}

{\bf Tau-plus.}
The following map can be uniquely extended to
an automorphism of $\HH\,$ 
(see \cite{C4},\cite{C101}):
\begin{align}
& \tau_+:\  X_b \mapsto X_b, \
\pi_r \mapsto q^{-\frac{(\om_r,\om_r)}{2}}X_r\pi_r,
\tau_+:\ T_0\mapsto X_0^{-1} (T_0+1),\ 
\label{taux}
\end{align}
This automorphism  fixes
$T_i\,(i\ge 1)$ 
and $q^{-1/(2m)},$ where the latter fractional power
of $q$  
must be added to the definition of $\HH\,.$  
\smallskip

\subsection{Polynomial representation}
The {\em Demazure operators}
are defined as follows:
\begin{align}
&T_i\  = \  (1-X_{\al_i})^{-1}(s_i-1),
\ 0\le i\le n;
\label{Demazx}
\end{align}
they obviously preserve $\Z[q][X_b, b\in P]$.
We note that only the formula for $T_0$ involves $q$:
\begin{align}
&T_0\  = (1-X_0)^{-1}(s_0-1),\hbox{\ where\ }\notag\\
&X_0=qX_\vth^{-1},\
s_0(X_b)\ =\ X_bX_{\vth}^{-(b,\vth)}
 q^{(b,\vth)},\
\al_0=[-\vth,1].
\end{align}

The map sending $ T_j$ to the corresponding operator from
(\ref{Demazx}), $X_b$ to the operator of multiplication
by $X_b$
(see (\ref{Xdex})), and 
$\pi_r (r\in O)$ to $\pi_r$ induces a
$ \Z_{q}$\~linear
homomorphism from $\HH\, $ to the algebra of
linear endomorphisms
of $\Z_q[X]$.
It will be called the
{\em polynomial representation};
the notation is
$$
\v\equal \ \Z_q[X_b, b\in P].
$$
It is faithful if $q$ is not a root of unity.

The polynomial representation
is actually the $\HH\,$\~module induced from 
the one-dimensional representation $T_i\mapsto 0,\,$ 
$\pi_r\mapsto 1, r\in O$
of the affine Nil-Hecke subalgebra 
$\h=\lan T_i,\pi_r\ran.$
\smallskip

{\bf Intertwiners.}
Let $T_i'\equal T_i+1$.
Given $\hw\in \hW$, 
the element $T'_{\hw}=\pi_rT'_{i_l}\cdots T'_{i_1}$ 
does not depend on the choice of the reduced
decomposition $\hw=\pi_r s_{i_l}\cdots s_{i_1}$.

We set $\Psi_{\hw}'=\tau_+(\pi_rT'_{i_l}\cdots T'_{i_1})$. 
Then $\Psi_i'\equal\tau_+(T'_{-\om_i})$ for
$i=1,\ldots,n$ are pairwise commutative and, importantly,
$W$\~invariant in the polynomial representation. 

Indeed, 
$\Psi_{b}'=\prod_{i=1}^n\, (\Psi_i')^{n_i}$ for
$P_-\ni b=-\sum n_i\,\om_i$. Provided that all $n_i>0$,
the {\em reduced} decomposition $b=b_-=w_0\pi_{b_+}$ holds
for the longest element $w_0\in W$ and $b_+=w_0(b)\in B_+$.
Thus $\Psi_b'$ is divisible on the left by $T_i'=T_i+1$
for any $i>0$ and therefore is divisible on the
left by the $W$\~symmetrizer. It results in 
the $W$\~invariance of $P_b$ for any $b\in P_-$.   
\smallskip

{\bf $q$\~Hermite polynomials.}
The most constructive way to define the
$q$\~Hermite polynomials is by using the intertwiners:
\begin{align}\label{barpform}
&P_b\ \equal\ q^{(b,b)/2}\Psi_b'(1) \,\for b\in B_-\,.
\end{align} 
Here $\Psi_i'$ can be replaced by their restrictions
to $\v^{\,W}$; they will become then 
pairwise commutative $W$\~invariant
difference operators. It is also the most suitable way
in this particular paper due to the direct connection with
the level-one Demazure characters.
\smallskip

It is important that the expansions of the $q$\~Hermite
polynomials in terms of $X_b$ and $q$ 
have only nonnegative (integral) coefficients.
This is the same for their nonsymmetric generalizations.
One can deduce the positivity from the interpretation of 
these polynomials via the Demazure characters 
or by direct using the DAHA-intertwiners. 

As $q\to 0$, the polynomials $P_{b}$ become
the classical finite dimensional Lie characters,
which can be seen, for instance, from (\ref{normppolsbar})
below.
\smallskip

\subsection{Inner products}
Let $\mu_\circ\equal \mu/\langle \mu \rangle$ for
\begin{align}
&\mu\ = \ \prod_{\al \in R_+}
\prod_{j=0}^\infty (1-X_\al q_\al^{j})
(1-X_\al^{-1}q_\al^{j+1})
\label{mubar}
\end{align}
and the constant term functional $\lan\,\cdot\,\ran$
(the coefficient of $\prod X_i^0$).
The following is well-known: 
\begin{align}
&\langle\mu\rangle\ =\ \prod_{i=1}^{n}
\prod_{j=1}^{\infty} \frac{1}
{1-q_i^j}\,, \where q_i=q^{\nu_i},\ 
\nu_i=\nu_{\al_i}=(\al_i,\al_i)/2.
\label{constermbar}
\end{align}

The polynomials $P_b (b\in P_-)$  can be uniquely determined from
the conditions:
\begin{align}\label{macd}
&P_b-\sum_{c\in W(b)}X_c\in 
\oplus_{c\succ b}\Z_{q} X_c \and
\lan P_b X_{c}\mu\ran = 0 
\for c\succ b,\\
&\{c\succ b\} \equal \{c\in W(b')\,\mid\, b'=b+a\in P_-
\hbox{\ for\ } 0\neq a\in Q_+\}. \notag 
\end{align}

For $b,c\in P_-\,$,
the norm formula reads as:
\begin{align}\label{normppolsbar}
&\lan\, P_{b} P_{c^\iota} \mu_\circ\,\ran\ 
=\ \de_{bc}\prod_{i=1}^n\prod_{j=1}^{-(\al_i^\vee,b)}
(1-q_i^{j})\for \ c^\iota\equal -w_0(c).
\end{align} 
\smallskip

{\bf Gauss-type inner products.}
Let us denote by $\xi$  the natural
projection $P\to P/Q$ and fix a nonempty subset 
$\varpi\subset P/Q$. If $\xi(b)\in \varpi$ then 
$\xi(c)\in \varpi$ for all monomials $X_c$ in $P_b$. 
We will use the symbol $\Join$ for
the whole $P/Q$.

We set:
\begin{align}
&\theta_\varpi(X)\ \equal\ 
\sum_{\xi(b)\in \varpi} q^{(b,b)/2}X_b\,, b\in P.
\label{gauser}
\end{align}

Due to \cite{C5} and \cite{ChW}, we obtain
the following formulas ($b,c\in P_-\,$):
\begin{align}
&\langle \theta_\varpi\mu \rangle\, =\, 1,\ 
\langle \theta_\varpi\mu_\circ\rangle\, =\,
\prod_{i=1}^n\prod_{ j=1}^{\infty}(1-q_i^j),
\label{mehtamul}\\
&\hbox{provided that }0\in \varpi \,\hbox{ and } 0 
\hbox { otherwise}, \notag\\
\label{pbgaussl}
&\langle P_b(X) P_{c^\iota}(X) 
\theta_\varpi\mu_\circ \rangle\,  =\,  
q^{\frac{(b-c)^2}{2}}
\langle \theta_\varpi\mu_\circ\rangle\\
&\hbox{for\ \,} 
\xi(c-b)\in \varpi \and 0\, \hbox {\ otherwise.} \notag
\end{align} 
Recall that $c^\iota=-w_0(c)$.
\smallskip

Compare with the ``free" formulas:
\begin{align}
&\langle \theta_\varpi\rangle\, =\, 1 \for
0\in \varpi \,\hbox{ and } 0 
\hbox { otherwise}, \notag\\
\label{pbgausslf}
&\langle X_b X_{-c} 
\theta_\varpi\rangle\,  =\,  
q^{\frac{(b-c)^2}{2}}\langle \theta_\varpi\rangle \for  b,c\in P
\\
&\hbox{when\ \ } 
\xi(c-b)\in \varpi \and 0\,\ \hbox {\ otherwise.} \notag
\end{align} 

Notice that $b,c$ here are in the whole $P$. Switching to
nonsymmetric $q$\~Hermite polynomials in (\ref{pbgaussl}),
provides better matching the free formulas. However, there
will be still some differences in the structure of 
these formulas.

\smallskip

\subsection{Main results}
The following theorem directly results from formulas
(\ref{normppolsbar},\ref{pbgaussl}).
We use the following  very basic fact from linear algebra and
functional analysis. 

Provided the convergence,
an arbitrary Laurent series $f(X)$ can be 
expressed as $\sum_b (\lr f, e'_b\rr/\lr e_b, e'_b\rr)e_b$ 
for any two bases $\{e_b\},\{e'_b\}$ in the space
of Laurent polynomials/series orthogonal to each other
with respect to a certain nondegenerate form $\lr\,,\,\rr$.

We consider expressions below as series in terms of 
nonnegative powers of $q$ (maybe fractional); analytically, 
one can assume that $|q|<1$. 
Recall that $q_i=q^{\nu_i}$
for $\nu_i=(\al_i,\al_i)/2$.

\begin{theorem} \label{MAINTHM}
Let us fix an arbitrary sequence of nonempty subsets 
$\boldsymbol{\varpi}=\{\varpi_i\in P/Q,\ 1\le i \le p\}$ and set
$\theta_i=\theta_{\varpi_i}$.
Then for the sequences $\mathbf b=\{b_i\in P_-, 1\le i\le p\}$\,:
\begin{align}\label{pggmain}
&\lan \mu \ran^p\,\prod_{i=1}^p \theta_i\ \,=\ \,
\sum_{\mathbf b}\frac{
q^{\bigl(\,b_1^2+(b_1-b_2)^2+\ldots+(b_{p-1}-b_p)^2\,\bigr)/2}}
{\prod_{i=1}^p\prod_{j=1}^n\prod_{k=1}^{-(\al_j^\vee,b_i)}
(1-q_j^{k})}\,
P_{b_p}(X),\\
&\hbox{subject to } \xi(b_1)\in \varpi_1,\, 
\xi(b_i-b_{i-1})\in \varpi_i \for 1<i\le p.\notag
\end{align}
In particular, for a fixed $b_p\in P_-$,
the corresponding sum on the right-hand side of (\ref{pggmain}) 
depends only 
on the unordered set $\{\varpi_i\}$. 
\sq 
\end{theorem}

The case of $p=1$ reads:
\begin{align}\label{pggmainone}
&\lan \mu \ran\,\theta_\varpi\ =\ 
\sum_{b\in P_-}\frac{
q^{b^2/2}}
{\prod_{j=1}^n\prod_{k=1}^{-(\al_j^\vee,b)}
(1-q_j^{k})}\,
P_{b}(X) \for \xi(b)\in \varpi.
\end{align}
Its iterations based on identities 
(\ref{normppolsbar}) and (\ref{pbgaussl})
readily provide the general (any $p$) formula.
\smallskip

The ``free" version of (\ref{pggmain}) based
on (\ref{pbgausslf}) is as follows\,:
\begin{align}\label{pggmainf}
&\prod_{i=1}^p \theta_i\ =  
\sum_{\mathbf b=\{\,b_1,\ldots,b_p\}}\,
q^{\bigl(b_1^2+(b_1-b_2)^2+\ldots+(b_{p-1}-b_p)^2\bigr)/2}
\,X_{b_p},\\
&\hbox{when\ \ } \xi(b_1)\in \varpi_1,\ 
\xi(b_i-b_{i-1})\in \varpi_i \hbox{\ \,for\ \,} 1<i\le p\,.\notag
\end{align}
The summation here is over all $b_1,\ldots,b_p\in P$
(not only anti-dominant).
\medskip

We will begin with the key application, which is
eliminating the $q$\~Hermite polynomials by
applying $\lan \,\cdot\, P_{c^\iota}\mu_\circ\,\ran$ 
in (\ref{pggmain}).
\medskip

{\bf Taking the constant term.}
\begin{corollary} \label{MAINCOR}
For a given $c\in P_-$, assuming that the summation variables 
$b_i$ are from $P_-$
and for an arbitrary sequence of nonempty subsets 
$\boldsymbol{\varpi}=\{\varpi_i\subset P/Q,\ 1\le i \le p\}$, 
which determine $\theta_i=\theta_{\varpi_i}$, 
\begin{align}\label{pggmcor}
&\Xi_{\boldsymbol{\varpi}}^{p,c}\equal\lan \mu \ran^p
\lan\,\prod_{i=1}^p{\theta_i}\, P_{c^\iota}\mu_\circ\ran
\!\!=\!\!\!\! 
\sum_{b_1,\ldots,b_{p-1}}\!
\frac{q^{\bigl(b_1^2+(b_1-b_2)^2+\ldots+
(b_{p-1}-c)^2\bigr)/2}}
{\prod_{i=1}^{p-1}\prod_{j=1}^n\prod_{k=1}^{-(\al_j^\vee,b_i)}
(1-q_j^{k})}, \\
&\xi(b_1)\in \varpi_1,\ 
\xi(b_i-b_{i-1})\in \varpi_i \hbox{ for } 
1<i<p,\ \xi(c-b_{p-1})\in \varpi_p\,.\notag
\end{align} 

Given an arbitrary collection $\boldsymbol{\varpi}$
of {\dfont atomic} $\varpi_{i}=\tilde{k}_i\equal \{k_i\}$
for $k_i\in O\cong P/Q$, there exists a unique $c=-\om_r$ 
for $r\in O$ (can be zero) such that using $P_{c^\iota}$ makes 
the left-hand side of (\ref{pggmcor}) nonzero, namely,
satisfying the congruence
$\,\om_r=\sum_{i=1}^{n} \om_{k_i}\mod Q\,$. Using the
notation $\Xi_{\boldsymbol{\varpi}}^{p,r}=
\Xi_{\boldsymbol{\varpi}}^{p,c}$ for $c=-\om_r$, one has:
\begin{align*}
\Xi_{\boldsymbol{\varpi}}^{p,r}
=q^{\frac{1}{2}\sum_{i=1}^{p-1} \om_{k_i}^2+
(\om_{k_1}+\ldots+\om_{k_{p-1}}-\om_r)^2/2}
\bigl(1+\sum_{j=1}^\infty C_j q^j\bigr) 
\hbox{ for } C_j\in \Z_+,
\end{align*}
where $\boldsymbol{\varpi}=\{\tilde{k}_i\}$ is atomic
and $r$ is defined as above.
\sq
\end{corollary}

Alternatively, the product 
$\prod_{i=1}^p{\theta_i}$ on the left hand-side of (\ref{pggmcor})
can be calculated using (\ref{pggmainf}) and $\mu_\circ$
can be replaced by its Laurent expansion, which can be
readily obtained directly using the formulas for the action 
of the lattice $Q$ on $\mu$ by translations.
Then we will arrive at a $q$\~series without 
the denominators. Multiplying it by $\lan \mu \ran^p\ $ and 
comparing with the right-hand side of (\ref{pggmcor}), we will obtain
Rogers-Ramanujan type identities.
One can also use here a developed machinery of finding formulas
for the products of theta functions in terms of the standard 
ones, more conceptually, the relations in the algebra of theta 
functions for the powers of an elliptic curve.

\comment{
For instance, the orbit-sums  
$$
q^{-lx^2/2}\,\sum_{\hw\in \hW}\,\hw(X_b\, q^{lx^2/2})
$$ 
can be taken, where we set $X_a=q^{x_a}$ and define 
$x^2$ accordingly. This will result
in formulas for the right-hand side of (\ref{pggmcor})
in terms of generalized $\eta$\~functions and similar
functions. More conceptually, we need to use here the
relations in the algebra of theta functions for the
elliptic curve and its powers.
The resulting formulas are an important
part of the theory of Rogers-Ramanujan identities 
of modular invariant type. We use this approach in the 
examples considered below.
}
\smallskip

\setcounter{equation}{0}
\section{Applications, modular invariance}
It is important to understand the actual range of 
our approach and find generalizations. Accordingly, we consider 
here not only the formulas that can be obtained by our method,
but also neighboring ones, those beyond our reach by now.

\subsection{Level-two formulas for B,C}\label{sec:lev2BC}
The general structure of
the multi-dimensional summations as in 
(\ref{pggmcor}) is of course not new. According to
\cite{An}, the existence of multi-dimensional 
generalizations of the classical Rogers-Ramanujan 
identities is common in the vast theory of the
Rogers-Ramanujan type identities. 

The case of $A_n$ and $\mathfrak{gl}$ will be examined below. 
The structure of our level-two formulas 
for the root systems $B_n,C_n$ is similar to that of
formulas (1.3),(1.8) from \cite{An},
but there are differences. 

Let us first  consider $B_n$.
In the notation from \cite{Bo}, our 
quadratic form is as follows:
$(b,b)=2\sum_{i=1}^n u_i^2$
for $b=\sum_{i=1}^n u_i \vep_i$. Also,
$$
(b,\al_j^\vee)=(b,(\vep_j-\vep_{j+1})/2)=u_j-u_{j+1}
\hbox{ for } j<n,\, 
(b,\al_n^\vee)=(b,\vep_n)=2u_n.
$$      
Recall that we normalize the inner product 
by the ``twisted" condition $(\al_{\sht},\al_{\sht})=2$ and
$\al_{\sht}^\vee=\al_{\sht}$; also $\,q_{\sht}=q,\ q_{\lng}=q^2$.
\smallskip

Let us take $\,c=0\,$ and  $\,\varpi_1=\{0\}=\varpi_2$. Setting
$\,v_i=-u_{n-i+1}$, one obtains: 
\begin{align}
&\lan \mu \ran^2\,
\lan\,\theta_1\theta_2\,\mu_\circ\,\ran
\! =\!\!\!\!\!\! 
\sum_{0\le v_1\le v_2\le \ldots \le v_n}
\frac{q^{\,2\,\sum_{i=1}^n v_i^2}}
{\prod_{k=1}^{2v_1}(1-q^k)\,
\prod_{i=1}^{n-1}\prod_{k=1}^{v_{i+1}-v_{i}}
(1-q^{2k})}.\notag
\end{align}

Compare with the summations in 
(1.3) and (1.8) from \cite{An}, where correspondingly
$d=1,2$\,:
\begin{align*}
\!\!\!\!\!\! 
\sum_{0\le v_1\le v_2\le \ldots \le v_n}
\frac{q^{\,\sum_{i=1}^n v_i^2}}
{\prod_{k=1}^{dv_1}(1-q^k)\,
\prod_{i=1}^{n-1}\prod_{k=1}^{v_{i+1}-v_{i}}
(1-q^{k})}.
\end{align*}

Let us now consider the root system $C_n$.
The quadratic form becomes
$(b,b)=\sum_{i=1}^n u_i^2$
for $b=\sum_{i=1}^n u_i \vep_i$. One has:
$$
(b,\al_j^\vee)=(b,(\vep_j-\vep_{j+1}))=u_i-u_{j+1} \for j<n,\ 
(b,\al_n^\vee)=(b,\vep_n)=u_n.
$$   
Taking $\,c=0,\ \varpi_1=\Join=\varpi_2\,$ and setting
$\,v_i=-u_{n-i+1}$ as above:
\begin{align}
&\lan \mu \ran^2\,
\lan\,\theta_1\theta_2\,\mu_\circ\,\ran
\! = \!\!\!\!\!\!  
\sum_{0\le v_1\le v_2\le \ldots \le v_n}
\frac{q^{\,\sum_{i=1}^n v_i^2}}
{\prod_{k=1}^{v_1}(1-q^{2k})\,
\prod_{i=1}^{n-1}\prod_{k=1}^{v_{i+1}-v_{i}}
(1-q^{k})}.\notag
\end{align}
The quadratic form here is as in \cite{An}, but there is
no match of the denominators.

It is instructional to consider the cases $B_1$ and $C_1$.
The corresponding $B_1$\~series is precisely the $A_1$\~series 
from (\ref{pggagact2mu}) for $\breve{\th}$ (for the lattice $Q$).
The $C_1$\~series becomes that for $\th$ there upon 
substitution $q^2\mapsto q$. So this is exactly as one 
can expect.

{\bf Product formulas.} Thanks to the Ole Warnaar, we can 
provide here the product expressions for the $q$\~series above. In the
$B$\~case, he used the Bailey pair $F(1)$:
\begin{align}\label{warnB}
\!\!\!\!\!\!  \sum_{0\le v_1\le v_2\le \ldots \le v_n}
&\frac{q^{\,2\,\sum_{i=1}^n v_i^2}}
{\prod_{k=1}^{2v_1}(1-q^k)\,
\prod_{i=1}^{n-1}\prod_{k=1}^{v_{i+1}-v_{i}}(1-q^{2k})}\\
=\ &\frac{(-q^{2n+1},-q^{2n+3}, q^{4n+4}\,;\, q^{4n+4})_\infty}
{(q^2\,;\,q^2)_\infty}\notag
\end{align}  
in the notation from (\ref{qn-not}) below. It generalizes
(\ref{pggagact2mu}), where $n=1$. 

The density of the terms in the numerator versus those in the
denominator can be readily calculated. It is
$\frac{6}{4n+6}=\frac{3}{2n+3}$, where 
$1-\frac{3}{2n+3}=\frac{2n}{2n+3}$, the density of the
missing terms in the numerators, is supposed to be 
connected with $c_{\hbox{\tiny eff}}$ for a dual
root system. It really coincides with that for $T_n$ in the 
introduction or in the table in Section \ref{sec:half-int} 
(without $\flat$). Actually this is a rigorous mathematical 
connection, which can be justified following \cite{Za,VZ}, but
we do not discuss this here in detail.

The $C_n$\~case is due to David Bressoud (we thank Ole Warnaar
for identifying it):
\begin{align}\label{warnC}
\!\!\!\!\!\!  \sum_{0\le v_1\le v_2\le \ldots \le v_n}
&\frac{q^{\sum_{i=1}^n v_i^2}}
{\prod_{k=1}^{v_1}(1-q^{2k})\,
\prod_{i=1}^{n-1}\prod_{k=1}^{v_{i+1}-v_{i}}
(1-q^{k})}\\
=\ &\frac{(q^{n+1},q^{n+1}, q^{2n+2}\,;\, q^{2n+2})_\infty}
{(q)_\infty}.\notag
\end{align}
 
Here the density of the terms in the numerator versus 
the denominator is obviously $3/(2n+2)$. Accordingly,
$1-\frac{3}{2n+2}=\frac{1}{2}\frac{2n-1}{n+1}$ is
$\frac{1}{2}$ times $c_{\hbox{\tiny eff}}$ for $B_n$
from Section \ref{sec:half-int}. 
The factor $\frac{1}{2}$ will disappear if the density
is recalculated to the base $q^2$ as in the previous formula; 
we will skip the details.
\smallskip

{\bf Warnaar's identities}.
Very close formulas to our $B,C$\~ones can be found
in \cite{War}. The formula right after (5.6) there
states that:
\begin{align}\label{warnaarB}
&\!\!\!\!\!\!  \sum_{0\le v_1\le v_2\le \ldots \le v_n}
\frac{q^{\sum_{i=1}^n v_i^2}}
{\prod_{k=1}^{2v_1}(1-q^k)\,
\prod_{i=1}^{n-1}\prod_{k=1}^{v_{i+1}-v_{i}}
(1-q^{2k})}= \prod_{j=1}^\infty 
\frac{1}{1-q^j}\,,\\
&\hbox{\ \ where\ \,} j\neq 2 \hbox{\ mod\,} (4)\and 
j\neq 0,\,\pm 4(n+1) \hbox{\ mod\,} (8n+12).\notag
\end{align}
Here and in (\ref{warnaarC}), the denominators exactly match 
our ones, but the quadratic forms are divided by $2$ versus 
our formulas. The next but one formula after (5.6) has 
the quadratic form extended by 
certain linear terms, which is analogous to our
using  $\,\varpi_1=\{1\}=\varpi_2$.

The density of the {\em missing\,} 
terms in the numerator of (\ref{warnaarB})
versus those in the denominator is
$\frac{3}{2} \frac{n+1}{2n+3}$, which is 
$\frac{3}{2}  (c_{\hbox{\tiny eff}}^\flat)_{C_n}$ 
for the (expected) formulas for $c_{\hbox{\tiny eff}}^\flat$ from
the table in Section \ref{sec:half-int}. We
will not comment on the factor $\frac{3}{2}$. 

A counterpart of our $C$\~case  
is his (5.14), which reads:
\begin{align}\label{warnaarC}
\!\!\!\!\!\!  \sum_{0\le v_1\le v_2\le \ldots \le v_n}
&\frac{q^{(\sum_{i=1}^n v_i^2)/2}}
{\prod_{k=1}^{v_1}(1-q^{2k})\,
\prod_{i=1}^{n-1}\prod_{k=1}^{v_{i+1}-v_{i}}
(1-q^{k})}\\
=\ \, &\frac{(q^{n/2+1/2}, q^{n/2+1}, q^{n+3/2}; q^{n+3/2})_\infty}
{(-q)_\infty (q^{1/2};q^{1/2})_\infty}\,.\notag
\end{align}

The density of the missing terms in the numerator
is $\frac{3}{4}\frac{2n-1}{2n+3}$, which is 
$\frac{3}{4} (c_{\hbox{\tiny eff}}^\flat)_{B_n}$,
where $\frac{3}{4}$ is  $\frac{3}{2}$ above divided by
$2$ due to the base $q^{1/2}$ in this case, similar to
the analysis of (\ref{warnC}).

There is a natural question.
Are identities (1.3) and (1.8) from \cite{An} and those
from \cite{War} somehow associated with root systems $B,C$? 

We are very thankful to
Ole Warnaar for letting us know about (\ref{warnaarB}) and 
(\ref{warnaarC}). We would like to mention
Section 5 from \cite{War} with interesting 
applications to certain Virasoro characters (which may be related
to our paper).
\smallskip

\subsection{Weyl algebra and Gaussian sums}
For an integer $N\ge 1$, we set 
$\ze=e^{\frac{2\pi\imath}{N}}$ and pick
$\ze^{1/(2m)}=e^{\frac{\pi\imath}{Nm}}$, 
where $(P,P)=\Z/m$ as above. Actually $\ze^{1/(2m)}$
can be any $(2m)$-th root of $\ze$, not necessarily a 
primitive $2mN$-th root (unless $2m$ divides $N$). 
Moreover, such nonprimitive
roots do appear and are necessary for the exact analysis of
the modularity; this is related to Theorem A
from \cite{KP}. In this section we stick to the choice above;
accordingly, the modular invariance condition we obtain will be
not always sharp. 
\smallskip

{\bf The extended Weyl algebra}. This algebra will be denoted by 
$\w_{N}$. It is generated over $\Q[\ze^{1/m}]$ by 
$U_a,V_b$ for $a,b\in P$
and  $w\in W$ subject to the relations,
$\,U_{a+b}=U_aU_b,\, V_{a+b}=V_aV_b$,
$$
wU_aw^{-1}=U_{w(a)},\ wV_aw^{-1}=V_{w(a)},\
U_aV_bU_a^{-1}V_b^{-1}=\ze^{(a,b)}.
$$

Setting $P[N]\equal P\cap NQ^\vee$, the quotient
\begin{align}\label{pnan}
\a_{N}\equal \Q[\ze^{1/m}][X_b,b\in P]/(X_c-1, c\in P[N])
\end{align}
has a natural structure of a $\w_{N}$\~module:
$$
U_b(X_a)=X_{a+b},\ V_b(X_a)=\ze^{-(a,b)}X_a \for a,b\in P.
$$
It can be canonically identified with the algebra
$$
\hbox{Funct}(P/P[N])=\oplus_{\bar a} \Q[\ze]\de_{\bar a}, \where
a\in P,\ \,\bar a =a \hbox{ mod } P[N].
$$ 
One has: $\de_{\bar a}\de_{\bar b}=\de_{\bar{a},\bar{b}}$
(the latter is the Kronecker delta),
$$
U_b(\de_{\bar a})=\ze^{(b,a)}\de_{\bar a},\ \,
V_b(\de_{\bar a})=\de_{\overline{a+b}},\ \,
\for
a,b\in P.
$$ 

Here and
below see Section 3.11.1, Section 3.10.4 from \cite{C101}
in the special case $k_{\lng}=1=k_{\sht}$, and also 
Lemma 3.11.3 there.
We will constantly use that $(\al^\vee,\al^\vee)=1$ for long 
$\al\in R$; the lattice $Q$ is even.
\smallskip

\begin{lemma}\label{pnlattice}
The coincidence $P[N]=NQ$ does not hold 
when and only when

\ \ \ (a) $N$ is odd for the systems $C_n(n\ge 2)$, $G_2$$(3|N)$,

\ \ \ (b) $N$ is even for the systems \,$B_n,\, C_n\,(n\ge 2)$,
$F_4$, \\
Moreover,
$P[N]=NQ^\vee$\, if\, 
$\nu_{\lng}\,|\, N$ or $N\in 2\Z_+$. Also,
$P[N]=NP$\, if\, $(N,\nu_{\lng})=1$ for 
$C_n, F_4, G_2, B_{n\in 2\Z}.$ 
Here $(\cdot\,,\cdot)=\,$gcd$(\cdot\,,\cdot)$ 
and by \,$a|b$\,, we mean
that $a$ divides $b$. \sq
\end{lemma}

The module $\a_N$ has a natural {\em projective} action of
$SL(2,\Z)$ that commutes with the action of $w\in W$ and
induces the standard action of $SL(2,\Z)$ on the generators
$\{U_a,V_b\}$.
Namely, for an arbitrary $\xi\in \C^*$, we set
\begin{align}\label{sl2triv}
&T_+=\left(\begin{array}{cc}1&  1 \\
                       0 & 1 \\
\end{array} \right)\,\rightsquigarrow\, \tau_+
(\de_{\bar a})\ =\ \frac{1}{\xi}\ze^{\frac{a^2}{2}}
\de_{\bar a} (a\in P),\\
&T_-=\left(\begin{array}{cc}1&  0 \\
                       1 & 1 \\
\end{array} \right)\,\rightsquigarrow\, \tau_-(X_a)\ =\ 
\xi\ze^{-\frac{a^2}{2}}X_a (a\in P).
\end{align}

To use these formulas without problems in $\a_N$, we need 
to check that $(b+Na)^2-b^2=2N(b,a)+N^2(a)^2$ is divisible 
by $2N$ for $a\in Q^\vee, Na\in P$. 
It obviously holds when $a\in Q$
(the latter is an even lattice) and for even $N$. 
In the $C$-$G$~case $(a)$
from Lemma \ref{pnlattice} (odd $N$, $3|N$ for $G_2$), 
we apply (\ref{sl2triv}) to the 
following set of generators for $P/P[N]=P/NP\, (C_n)$ and 
$P/P[N]=Q/NQ^\vee\, (G_2)$:
\begin{align}\label{CnGgauss}
&P\ni b\ =\ \sum_{i=1}^n c_i\om_i,\ \, 0\le c_i <N 
\hbox{\ \, for \ } C_n (n\ge 2),\\
&Q\ni b=c_1\al_i+c_2\al_2,\ 0\le c_1 <N,\, 0\le c_2 <N/3
\for G_2.\notag
\end{align}

Let $\si\equal\tau_+ \tau_-^{-1}\tau_+$. Then one checks that
$$
\tau_+ \tau_-^{-1}\tau_+=\tau_-^{-1} \tau_+\tau_-^{-1},\ \, 
\si\tau_+^{\pm 1}\si^{-1}=\tau_-^{\mp 1},\ \,(\si\tau_+^{-1})^3=\si^2.
$$
These are the relations of the projective $PSL(2,\Z)$
due to Steinberg (the last two formally follow from the first). 
Explicitly,
\begin{align*}
&\si:\ 
\de_{\bar a}\mapsto \frac{\ga}{\xi^3} X_a, \,
\,X_a\mapsto \frac{\ga}{\xi^3}\de_{-\bar a}
\for \ga\equal\frac{\sum_{b\in P/P[N]}\ze^{\frac{b^2}{2}}}
{\hbox{\tiny $\sqrt{[P:P[N]]}$}}\,.
\end{align*}

In terms of $\de_{\bar a}$, the formula for $\si$ becomes
as follows:
\begin{align}\label{simatrix}
&S=\left(\begin{array}{cc}0& 1 \\
                      -1 & 0 \\
\end{array} \right)\rightsquigarrow\,
\si(\de_{\bar a})= 
\frac{\ga}{\xi^3}\sum_{b}\ze^{(a,b)}\de_{\bar b}
\for a,b\in P/P[N].
\end{align}

{\bf Gaussian sums.}
We will not discuss here the calculation of the Gaussian
sums and the corresponding $\ga$, which are always roots
of unity. 
For $A_1$, one obtains that $\ga=e^{\frac{\pi \imath}{4}}$,
the most involved instance of the celebrated four Gauss formulas.
We note that the formulas for Gaussian sums can be deduced
from the $q$\~Mehta-Macdonald identities;  see \cite{C101}, 
Sections 2.1 and 3.10. Let us provide the list
of $\ga$ for all root systems. They influence the conductors
of the modular functions under consideration if we need to know
the exact invariance properties, not up to a character, 
which is addressed in Theorem \ref{MODULARCOR} below.

\begin{proposition}\label{LEMALT}
(i) Setting $\varrho\equal e^{\frac{\pi \imath}{4}}$,
\begin{align}\label{gaussroots}
&\ga=\varrho^n \for A_n, D_n, E_{n=6,7,8},\\
F_4: \ga=(-1)^{N-1}&,\ \,
G_2: \ga\,=\,\imath,\, 1,-1 \for N\hbox{ mod }3=0,1,2,\notag\\
C_n (n\ge 2)&: \ga=\varrho^{n-\psi}, 
\where \psi=0 \hbox{\ \ unless: }\notag\\
&\psi=n \for N=1\hbox{ mod } 4 \and \notag\\ 
&\psi=1 \for \{N=3\hbox{ mod } 4 \ \&\, \hbox{ odd } n\},\notag\\
B_n (n\ge 3)&: \ga=\varrho^{n-\psi}, 
\where \psi=0 \hbox{\ \ unless }\notag\\
&\psi=4 \for N=1\hbox{ mod } 2.\notag
\end{align}

(ii) In the notation above, let us take $\xi^3=\pm \imath\ga$,
for instance, $\xi=\pm\varrho=\pm e^{\frac{\pi \imath}{4}}$ 
for $A_1$.
Then the operator $\si^2$ becomes the inversion
in Funct$(P/P[N])$, namely, 
$$
\si^2:\ \de_{\bar a}\mapsto -\de_{-\bar a},\ X_a\mapsto -X_a^{-1}.
$$
Accordingly, formulas (\ref{sl2triv}) and (\ref{simatrix})
induce the action of the group $PSL(2,\Z)$ in the 
space $\b_{N}\equal
\{f\in \a_N\, \mid\, w(f)=\sgn(w)z\}$ of the skew-symmetric 
elements of $\a_{N}$. \sq
\end{proposition}
\smallskip

In the simplest case of $p=1$, i.e. for $N=h+1$ the formulas
(\ref{gaussroots}) for the simply-laced root systems
follow from general facts about the Fourier transforms for
arbitrary even lattices. See, e.g. 
paper \cite{VW} for the case of $\mathfrak{sl}_n$ and \cite{Wu}
(esp,., Theorem 3.1 and formula (3.9) there).
Since, we verified that $\ga$ does not depend on $N$ for
$A$-$D$-$E$, this is sufficient to catch its value. 
\smallskip

\subsection{The action of SL(2,Z)}
By modular invariant functions
with respect to a congruence subgroup
$\Ga\subset SL(2,\Z)$ (actually its image
in $PSL(2,\Z)$), we mean the modular functions (of weight zero)
for a certain finite character of $\Ga$. It is up to 
roots of unity, if an individual $g\in SL(2,\Z)$ is considered.
The action of 
{\tiny $\left(\begin{array}{cc}
a & b \\
c & d \\
\end{array}\right)$}
$\in SL(2,\Z)$ is standard: $z\mapsto \frac{az+b}{cz+d}$\, for\, 
$q=e^{2\pi \imath z}$.

\begin{theorem}\label{MODULARCOR}
Given level $p$ as in the theorem, we set 
$N=p+h$ for the Coxeter number $h=(\rho,\vth)+1$ of the root system
$R$ and denote by $\Ga^{(N)}$ the kernel of the map 
$PSL(2,\Z)\to $Aut$(\b_N)/\C^*$:  
\begin{align}\label{projinb}
&T_+\mapsto \tau_+,\ T_-\mapsto \tau_-,\ 
S\mapsto \si \for \xi^3=\pm\imath\ga.
\end{align}
For instance, the elements $\tau_{\pm}^M$ for
$M=2Nm$ belong to $\Ga^{(N)}$.

Let us fix a minuscule $c=-\om_r\in P_-$ or take $c=0 (r=0)$
and consider all collections
$\boldsymbol{\varpi}=\{\varpi_i\in P/P[1]=P/(P\cap Q^\vee)\}$. 
Upon multiplication by
a proper common (depending on $p,R,c$) fractional power  
$\exp(2\pi \imath z\,e/f)$ of $\,q=\exp(2\pi \imath z)\,$ 
for $\,e,f\in \Z\, (f>0)\,$,
the $q$\~series 
\begin{align}\label{prodthmu}
\Xi_{\boldsymbol{\varpi}}^{p,r}\equal
\lan \prod_{i=1}^p{\theta_i} P_{c^\iota} \mu_\circ\ran\,
\lan \mu \ran^p
\end{align}
from Corollary \ref{MAINCOR} become modular $\Ga^{(N)}$\~invariant 
functions, which are also modular invariant with 
respect to $T_+^{m}$ (up to a character depending on $c$ and,
accordingly, on $e/f$). \sq
\end{theorem}
\smallskip
 
{\bf On the justification.} It essentially goes as follows. 
We define the {\em skew-symmetric Looijenga space} $\l^{\bar{b}}_N$
of level $N=p+h$  for $\bar b\in P/P[N]$
as the linear span of the orbit-sums over $Q^\vee\rtimes W$,
\begin{align*}
&\chi_b=\sum_{aw\in Q^\vee\rtimes W}\!\!\!\! 
(-1)^{l(aw)}\,aw(X_b q^{N\frac{x^2}{2}}) 
\hbox{\, provided that \,} Na\in P,\\
&\hbox{where\, formally\ \ } 
(aw)(q^{N\frac{x^2}{2}})\,=\,q^{N\frac{(x+a)^2}{2}}\,=\,
q^{N\frac{a^2}{2}}X_{Na}\,q^{N\frac{x^2}{2}}.
\end{align*}
Notice that $a\in Q^\vee$ here; $b$ is any pullback 
of $\bar b$ to $b\in P$. Let $\tilde{\chi}_b$
be $\chi_b$ divided by the coefficient of $X_{c}$ in 
$\chi_b$ (if it is nonzero) for $\,c\,$
from (\ref{prodthmu}).
 
The Kac-Moody characters of (twisted) irreducible 
level-$p$ integrable modules constitute a basis in 
$\l^{\bar{b}}_N$ upon their multiplication by $\mu$ (and the 
standard eta-type factors in terms of $q$).
Concerning the functional equations for the theta functions
associated with root systems and  the action of $PSL(2,\Z)$ on the
Kac-Moody characters and string functions  
see, e.g, Proposition 3.4 and Theorem A from \cite{KP} and 
\cite{Kac}.


Then we decompose 
$\lan\mu\ran^p(\prod_{i=1}^p{\theta_i}) P_{c^\iota} \mu_\circ\,$ in 
terms of such $\tilde{\chi}_b$. The modular 
invariance of the latter under the action of 
$g\in PSL(2,\Z)$ (up to proportionality) implies that for the 
coefficients in this decomposition. We use that 
the coefficients of $X_c$ were made $1$ or zero in 
$\tilde{\chi}_b$. The modular action of $PSL(2,\Z)$ 
in the skew-symmetric Looijenga space is
described in Proposition \ref{LEMALT} for the skew-symmetric
$W$\~orbit-sums there. It is upon proper normalization, which
results in a common factor $q^{e/f}$ for 
$\Xi_{\boldsymbol{\varpi}}^{p,r}$ (given $c$, for all $\varpi$).

Since our series are in terms of integral powers of
$q^{1/m}$, the modular invariance of $\Xi_{\boldsymbol{\varpi}}^{p,r}$
with respect to $T_+^{m}$ is granted (up to roots of unity). 

We note that the congruence subgroups 
of all such $g$ can be greater than $\Ga^{(N)}$ (extended by
$T_+^{m}$) for certain choices of $R,p$ and $\mathbf {\varpi}$. 
Recall that we need to consider the action of $PSL(2,\Z)$ only in the 
subspace $\b_N$. Also, picking 
a primitive root $\ze^{1/(2m)}$ is sufficient but not
always necessary.   Presumably, taking the (nonprimitive) 
root $\ze^{1/(2m)}=e^{\frac{2\pi\imath k}{Nm}}$ 
for odd $N=2k-1$ is always sufficient, but we did not check 
all details; cf. Theorem A from \cite{KP}. 
If this is true, then the special treatment of odd $N$ in 
the cases of $C_n$ and $G_2$ in (\ref{CnGgauss}) 
will be unnecessary.
\smallskip

Let us give the simplest example. When $p=1$, the set $\varpi$ can be
only $\tilde{0}\equal \{0\}$ in the absence of $P_c$ for minuscule
$c$; also $N=h+1$ is relatively prime with $m$ in this case. 
The subspace of $W$\~anti-invariant elements in Funct$(P/NP)$ 
is sufficient here instead of Funct$(P/P[N])$. Since this space  
is one-dimensional, the modular invariance must hold for 
the whole $PSL(2,\Z)$. This analysis is of course not necessary
because we know that $\lan \th_{\tilde{0}} \mu\ran=1$.

We also mention that our series for the greatest 
$\boldsymbol{\varpi}=\{\Join,\ldots,\Join\}$ are generally modular 
invariant for smaller powers of $T_-$ than those for generic 
$\boldsymbol{\varpi}$. 
\medskip

\subsection{The A-case}
Let us discuss Corollary \ref{MAINCOR}
for $A_{n}$. Then $\xi:\,P\to P/Q=\Z_{n+1}$ and
for $\{b_i\}\subset (P)_-\,$,
\begin{align}\label{pggmdual1}
&\frac{\lan\prod_{i=1}^{p}\theta_i\,P_{c^\iota}\,\mu_\circ\ran}
{(\prod_{j=1}^\infty (1-q^j))^{pn}}=\!\!
\sum_{b_1,\ldots,b_{s}}\!
\frac{q^{\bigl((b_1)^2+(b_1-b_2)^2+\ldots+
(b_{p-1}-c)^2\bigr)/2}}
{\prod_{i=1}^{p-1}\prod_{j=1}^{r}\prod_{k=1}^{-(\al_j,b_i)}
(1-q^{k})},\\
\xi(b_1)\in \varpi_1&,\ 
\xi(b_i-b_{i-1})\in \varpi_i,\, 1<i\le p-1,\ 
\xi(c-b_{p-1})\in \varpi_p.\notag
\end{align}
\medskip
Here $c=-\om_r\, (r\in O)\,,\,$ $\om_i^\iota=\om_{n-i+1}$.
This correction is necessary to
make the summation on the right-hand side nonzero;
given $\boldsymbol{\varpi}$, such $c$
is determined uniquely.
\smallskip

{\bf Level-rank duality.}
When $\varpi_i=\Join$ for all $i$, we come to the simplest
case of the level-rank duality. Let us divide
(\ref{pggmdual1}) by $(q)_\infty^{p-1}$
for  $(q)_\infty\equal \prod_{j=1}^\infty (1-q^j)$.
Then the right-hand side {\em coincides} with the 
Rogers-Ramanujan type expression for the $q$\~character of 
the module $L((n+1)\La_0)$ of $\hat{\mathfrak{sl}}_{p}$.
See formula (5.7) from \cite{Geo} and Theorem 7 from 
\cite{Pr} (and references therein). The origin of this
approach is due to \cite{SF}. 

We note that the division by $(q)_\infty^{p-1}$ makes
perfect sense, since the  resulting power
$(q)_\infty^{p(n+1)-1}$ on the left-hand side of  
(\ref{pggmdual1}) will become level-rank symmetric,
i.e. invariant with respect to $p\leftrightarrow n+1$.

We will not discuss in this work  the level-rank duality
beyond this observation and the case $n=1,p=3$ considered 
in Section \ref{sec:casep=3}.

\medskip

{\bf Switching to GL.}
Let us discuss a 
$\mathfrak{gl}_{n+1}$\~version
of formula (\ref{pggmdual1}). The notations are as follows:
\begin{align*}
&\R^{n+1}\ =\ \oplus_{i=1}^{n+1}\R\vep_i,\,\  
(b,b)=\sum_{i=1}^{n+1} (u_i)^2 \for 
b=\sum_{i=1}^{n+1} u_i \vep_i,\\
&\al_i=\vep_i-\vep_{i+1}\,(i\le n),\  \om_i=\vep_1+\cdots+\vep_i
\,(i\le n+1),\ \om_0=0, \\
&P=\oplus_{i=1}^{n+1}\Z\vep_i,\ Q=\oplus_{i=1}^{n}\Z\al_i,\ 
P_+=\{b\in P\,\mid\, u_i\ge u_{i+1}\},\ P_-=-P_+.
\end{align*}
We define $\xi(b)\equal \sum_{i=1}^{n+1} u_i \!\mod\! 
(n+1)\!\in \Z_{n+1}$ 
for $b\in P$.  Then 
\begin{align}\label{mubargl}
&\mu\ = \ \prod_{1\le i<j \le n+1}\,
\prod_{k=0}^\infty (1-(Z_i/Z_j) q_\al^{k})\,
(1-(Z_j/Z_i)q_\al^{k+1}),\\
&\theta_\varpi(X)\ \equal\ 
\sum_{\xi(b)\in \varpi} q^{(b,b)/2}Z_b\,,\ b\in P
\for \varpi\subset \Z_{n+1}.
\label{gausergl}
\end{align}
Here $b=\sum_{i=1}^{n+1} u_i \vep_i,\ 
Z_b=\prod_{i=1}^{n+1} Z_i^{u_i}\,$; 
the constant term will be defined with respect to
the variables $\{Z_i\}$. 

We pick $\theta_i=\theta_{\varpi_i}$ 
for a given sequence of subsets
$\boldsymbol{\varpi}=$
$\{\varpi_i\subset \Z_{n+1}\}$, where $i=1,2,\ldots,p$.
Let $c=-\om_r$ for $0\le r\le n$. Then $P_{c^\iota}$ is 
the $W$\~orbit-sum of $(Z_1\cdots Z_{n-r+1})^{-1}$ for $r>0$ and $1$
for $r=0$.
 
Taking the summation elements 
$\{b_i\}$ from $P_-$, we claim that
\begin{align}\label{pggmduagl}
&\frac{\lan\prod_{i=1}^{p}\theta_i\, P_{c^\iota}\,\mu_\circ\ran}
{\prod_{j=1}^\infty (1-q^j)^{pn}}\!=\!\!\!\!
\sum_{b_1,\ldots,b_{p-1}}\!\!\!\!
\frac{q^{\frac{1}{2}\bigl(b_1^2+(b_1-b_2)^2+\ldots+(b_{p-2}-b_{p-1})^2+
(b_{p-1}-c)^2\bigr)}}
{\prod_{i=1}^{p-1}\prod_{j=1}^{n}\prod_{k=1}^{-(\al_j,b_i)}
(1-q^{k})},\\
&\xi(b_1)\in \varpi_1,\ \,
\xi(b_i-b_{i-1})\in \varpi_i,\, 1<i<p,\ \,
\xi(c-b_{p-1})\in \varpi_p
\,.\notag
\end{align}
The left-hand side and therefore the right-hand side
do not depend on the order of $\varpi_i$ in the collection
$\boldsymbol{\varpi}.$ 
\smallskip

Let us establish a relation to  $A_{n}$.
Recall that the fundamental weights for  $A_{n}$
are connected with those for $\mathfrak{gl}_n$ as follows:
$$
\om_i'=\sum_{j=1}^i \vep_i-\frac{i}{n+1}\bar{\om},\ \,\bar{\om}=
\om_{n+1}=\vep_1+\ldots+\vep_{n+1},\ 1\le i \le n.
$$

Upon passage to the fundamental weights of $A_{n}$,
the connection between formulas (\ref{pggmduagl}) and 
(\ref{pggmdual1}) 
involves a nontrivial $\eta$\~type factor, which we
are going to calculate.

We use that the square $(b'+\bar{u}\bar{\om})^2$ defined for
$\mathfrak{gl}_{n+1}$ is $(b')^2+(n+1)\bar{u}^2$ for $b'$ 
expressed in terms of $\om_i'\, (i\le n)$. Here
$(n+1)\bar u\hbox{ mod} (n+1)$ equals $\xi(b)$ for $b\in P$ and 
coincides with $\xi(b')$ defined for $A_{n}$.

Let us consider only the {\em atomic} sets
$\tilde{k}=\{k\}\subset \Z_{n+1}$ for $k\in \Z_{n+1}$.
Then the collections 
$\,\boldsymbol{\varpi}\,$ $\,=\,\{\varpi_i,1\le i\le p\}\,$
are given by the decompositions $\,p=\la_0+\cdots+\la_{n}\,$,
where we treat $\la_k$ as the multiplicity of $\tilde{k}$ in 
$\boldsymbol{\varpi}$. The formula (\ref{pggmduagl}) can be
presented as follows:

\begin{align}\label{pgggltoa}
&\lan\prod_{i=1}^{p}\theta_i\,P_{c^\iota}\,\mu_\circ\ran\,
\lan \mu \ran^{p}\ =\ \prod_{j=0}^{n}\,\ 
(\sum_{s\in j/(n+1)+\Z}q^{s^2/2}\,)^{\la_j}\\
\times\sum_{b'_1,\ldots,b'_{p-1}}
&\frac{q^{\bigl((b'_1)^2+(b'_1-b'_2)^2+\ldots+(b'_{p-2}-b'_{p-1})^2+
(b'_{p-1}-c')^2\bigr)/2}}
{\prod_{i=1}^{p-1}\prod_{j=1}^{n}\prod_{k=1}^{-(\al_j,b'_i)}
(1-q^{k})}\,, \where\notag \\
\xi(b'_1)\in &\varpi_1,\ \,
\xi(b'_i-b'_{i-1})\in \varpi_i,\, 1<i<p,\ \,
\xi(c'-b'_{p-1})\in \varpi_p\,,
\,\notag
\end{align}
and the summation is in terms of arbitrary  $\{b_i'\}$
from the lattice $P_-'$ defined for $A_n$. Thus the sum on the
right-hand side is exactly as it was in the case of $A_{n}$.
\smallskip

\subsection{Rank 2 level 2} 
Let us consider the first nontrivial example in the case of 
$A_2$. See \cite{KKMM} for some coset aspects of this case.
We take $p=2\,$ and $c=0\,$; then there
can be only two admissible atomic collections 
$\boldsymbol{\varpi}=\{\varpi_1,\varpi_2\}$
in  $P/Q=\Z_3$, namely,
$$
\circledcirc=\{\tilde{0},\tilde{0}\}, \ \,
\circledast=\{\tilde{1},\tilde{2}\},\ \ 
\tilde{i}\equal\{i\}\subset \Z_3.
$$ 
The functions $\Xi_{\boldsymbol{\varpi}}^{2,0}$ for any other 
collections $\boldsymbol{\varpi}$ can be linearly expressed in 
terms of these two. For instance,
$$
\Xi_{\tilde 2,\tilde 1}^{2,0}\,=\,\Xi_{\tilde 1,\tilde 2}^{2,0},\ \
\Xi_{\Join}
\equal \Xi_{\Join,\Join}^{2,0}\,=\,
\Xi_{\tilde 0,\tilde 0}^{2,0}+2\,\Xi_{\tilde 1,\tilde 2}^{2,0}\,.
$$
Thus it suffices to consider the following:
\begin{align}\label{pggmdual20}
&\Xi_{\circledcirc}\!=\!
\lan\prod_{i=1}^{2}\theta_{\tilde{0}}^2\,\mu_\circ\ran\,
\lan \mu \ran^{2}\!=\!
\sum_{b}
\frac{q^{b^2}}
{\prod_{j=1}^{2}\prod_{k=1}^{-(\al_j,b)}
(1-q^{k})},\ b \in P_-\cap Q,\\
\label{pggmdual21}
&\Xi_{\circledast}\!=\!
\lan\prod_{i=1}^{2}\theta_{\tilde{1}}\theta_{\tilde{2}}\,
\mu_\circ\ran\,
\lan \mu \ran^{2}\!=\!
\sum_{b\in P_-}
\frac{q^{b^2}}
{\prod_{j=1}^{2}\prod_{k=1}^{-(\al_j,b)}
(1-q^{k})}\,,\ \xi(b)=1.
\end{align}

The sum on the right-hand side of (\ref{pggmdual20})
plus that in (\ref{pggmdual21}) times $2$, which
is $\Xi_{\Join}$,
can be found in \cite{Za} (Table 2, pg. 47)
and in \cite{VZ}, Table 1. It is
for the matrix $A=$
{\tiny $\left(\begin{array}{cc}
4/3 & 2/3 \\
2/3 & 4/3 \\
\end{array}\right)$} there and for $B=(0,0)^{tr}$.
Accordingly, $b^2=\frac{2}{3}(v_1^2+v_1v_2+v_2^2)$
for $b=v_1\om_1+v_2\om_2$.
It was expected in \cite{VZ}, based on 
computer calculations, that 
\begin{align}\label{pggmdual22}
&q^{-1/30}\,\sum_{b\in P_-}
\frac{q^{b^2}}
{\prod_{j=1}^{2}\prod_{k=1}^{-(\al_j,b)}
(1-q^{k})}\\
&=\,
\frac{1}{\eta(z)}
\sum_{n\in \Z}(-1)^n\bigl(
2q^{\frac{15}{2}(n+\frac{3}{10})^2}+
q^{\frac{15}{2}(n+\frac{1}{30})^2}-
q^{\frac{15}{2}(n+\frac{11}{30})^2}\bigr).\notag
\end{align}
However this formula was not finally confirmed there. 
Using $\Xi_{\Join}$,
we obtain that the left hand-side of
(\ref{pggmdual22}) is a 
modular function. Then one needs to compare only 
few terms in the $q$\~expansions to establish their
coincidence. It is a 
simple computer verification (which we performed).
\smallskip

Let us see what 
Corollary \ref{MODULARCOR} gives in this case. It states 
that all $\Xi$ are modular for at least $\Ga(30)$
up to a fractional power of $q$ (which is $-1/30$).
Exact formulas for these fractional powers of $q$ are
not discussed in our paper; cf. \cite{Za},\cite{VZ}.
Also, the modular invariance of 
$$
{}^\sharp \Xi\equal q^{-1/30}\Xi,
\where  \Xi \hbox{ \ equals\ }
\Xi_{\circledcirc},\  \Xi_{\circledast},\ 
\Xi_{\Join},
$$
holds at least with respect to
$T_+^m$ (which is obvious) and for $T_{-}^M$ for
$M=2Nm$. 
Recall that $T_-=${\tiny$\left(\begin{array}{cc}
1& 0 \\
1&1 \\
\end{array}\right)$.}

Since 
$\,N=p+h=5,\, m=3\,$
in this example, the corollary gives that at
least $M=30$ will be sufficient.
However we need to consider only
$\b_N$ (not the whole $\a_N$)\,
in this corollary and also $\,\rho\in Q\,$ for 
the root system $A_2$.
Therefore $M=15$ is actually sufficient here. 
This matches our expectation that $Nm$ instead of 
$2Nm$ gives the right (at least more exact) 
estimate for $M$ in the case of odd $N$.

We note that $M=5$ is sharp here for the $T_-^M$\~invariance 
of ${}^\sharp\Xi_{\Join}=\,q^{-1/30}\,\Xi_{\Join}\,$.
However,  
it is exactly $M=15$ for ${}^\sharp\Xi_{\circledcirc}$ and 
${}^\sharp\Xi_{\circledast}$ (as we computed directly).
Recall that by ${}^\sharp$, we mean the multiplication
by $q^{-1/30}$.
Moreover,
\begin{align}\label{5action2}
(T_-^{-5})\,{}^\sharp\!\overrightarrow{\Xi}\ =\ 
\left(\begin{array}{cc}
    \frac{1}{2}+\frac{\imath}{\sqrt{12}} & 
    -\frac{4\imath}{\sqrt{12}}\\
    -\frac{2\imath}{\sqrt{12}} &\frac{1}{2}-
   \frac{\imath}{\sqrt{12}} \\
  \end{array}\right)
{}^\sharp\!\overrightarrow{\Xi} \for
{}^\sharp\!\overrightarrow{\Xi}=
\left(\begin{array}{c}
    {}^\sharp\Xi_{\circledcirc}\\
    {}^\sharp\Xi_{\circledast}\\
\end{array}\right)\,.
\end{align}
In particular,
\begin{align}\label{2lev2tot}
(T_-^{-5})\, {}^\sharp\Xi_{\Join}\ =\ 
(T_-^{-5})\, ({}^\sharp\Xi_{\circledcirc}+
2\,{}^\sharp\Xi_{\circledast})\ =\ 
(\frac{1}{2}-\frac{3\imath}{\sqrt{12}})\, 
{}^\sharp\Xi_{\Join}\,,
\end{align}
where $\frac{1}{2}-\frac{3\imath}{\sqrt{12}}=e^{-2\pi\imath/6}$.
\smallskip

Here we use the expansions
of $\Xi_{\circledcirc}$ and
$\Xi_{\circledast}$:

\begin{align}\label{pggmdual220}
&{}^\sharp\Xi_{\circledcirc}\ =\ \frac{1}{\eta(z)}
\sum_{n\in \Z}(-1)^n\bigl(
q^{\frac{15}{2}(n+\frac{1}{30})^2}-
q^{\frac{15}{2}(n+\frac{11}{30})^2}\bigr),\\
\label{pggmdual221}
&{}^\sharp\Xi_{\circledast}\ =\ \frac{1}{\eta(z)}
\sum_{n\in \Z}(-1)^n\bigl(
q^{\frac{15}{2}(n+\frac{3}{10})^2}\bigr).
\end{align}
Using Corollary \ref{MAINCOR} or directly, we observe that
$$
 \Xi_{\circledcirc,\circledast}=
q^{\om_{k}^2}
\,\bigl(1+\sum_{j=1}^\infty C_j q^j\bigr), 
$$
where $k=0,1$  correspondingly. This can be readily
checked for the right-hand sides in (\ref{pggmdual220}
\ref{pggmdual221}) multiplied by $q^{1/30}$;
the smallest powers of $q$ become respectively 
$\frac{1}{30}-\frac{1}{24}+\frac{1}{120}=0$ and
$\frac{1}{30}-\frac{1}{24}+\frac{27}{40}=\frac{2}{3}$.

We note that the formulas in Table 1 from
\cite{VZ} with nonzero matrix $B$ (for the same $A$
as above) correspond to taking nonzero minuscule $c$
in our construction, i.e. to
$\Xi^{2,r}_{\Join}$ for $r=1,2$. 
\smallskip

{\bf Level-rank duality.}
Following Section \ref{sec:casep=3} below
(in notation from \cite{Geo}),
the level-rank duality predicts that 
$$
{}^\sharp\Xi_{\circledcirc}/(q^{-1/24}\eta(z))
\and {}^\sharp\Xi_{\circledast}/(q^{-1/24}\eta(z))
$$
must coincide with the following string functions for 
$\hat{\mathfrak{sl}}_2$ of level $3$\,:
\begin{align*}
q^{-1/120}\frac{\Xi_{\circledcirc}}{\eta(z)}\ =\ 
c_0^{3\hat{\La}_0},\ 
q^{-1/120}\frac{\Xi_{\circledast}}{\eta(z)}\ =\ 
c_{\al_1}^{3\hat{\La}_0}.
\end{align*}
The level-rank duality for any $\,R,p$\, is
not discussed in this work.

Indeed, we checked that 
\begin{align}\label{pggmdual2200}
&\ \ \ \ \ q^{-1/120}\,(\Xi_{\circledcirc}-
\Xi_{\circledast})\,/\,\eta(z)\\ 
&=\ q^{-1/120}\sum_{n\in \Z}(-1)^n\bigl(
q^{\frac{15}{2}(n+\frac{1}{30})^2}-
q^{\frac{15}{2}(n+\frac{11}{30})^2}
-q^{\frac{15}{2}(n+\frac{3}{10})^2}\bigr)\notag
\\
&=\ \prod_{k=1}^\infty (1-q^{k/3}) \where
k\neq \pm 1 \hbox{\, mod\,\,}(5).\notag
\end{align}
The latter product exactly coincides with
$c_0^{3\hat{\La}_0}-c_{\al_1}^{3\hat{\La}_0}\ =\ 
c_{30}^{30}-c_{12}^{30}$\,,
where the notation and the formula for the last 
difference can be found at page 220 in \cite{KP}.  


\medskip

\setcounter{equation}{0}
\section{The rank-one case}
Here $\al=\al_i=\vth$, $s=s_1$, $\om=\om_1=\rho$; so
$\al=2\om$ and the standard invariant form reads as
$(n\om,m\om)=nm/2$.
Also, we set
$$
X=X_\om=q^x,\ X(q^{n\om})=q^{n/2},\ 
\Ga(F(X))=\om^{-1}(F(X))=F(q^{1/2} X), 
$$
i.e. $x(n\om)=n/2$, $\Ga(x)=x+1/2$.
Also, $\pi\equal s\Ga: X\mapsto q^{1/2}X^{-1}$;


\subsection{Basic functions}
The $q$\~Hermite polynomials will be denoted by
$P_n \equal P_{-n\om}\, (n\in \Z_+)\,$ in this section.
For instance, $P_0=1$,\ $P_1=X+X^{-1},$
\begin{align}
P_2=&X^2+X^{-2}+1+q,\ 
P_3=X^3+X^{-3}+\frac{1-q^3}{1-q}(X+X^{-1}),\notag\\
&P_4=X^4+X^{-4}+\frac{1-q^4}{1-q}(X^2+X^{-2})
+\frac{(1-q^4)(1-q^3)}{(1-q)(1-q^2)}.\label{overpsmall}
\end{align}
The general formula is classical. For the monomial 
symmetric functions $M_0=1$, $M_n=X^n+X^{-n}$ $(n>1)$, 
\begin{align}
&P_n=M_n+\sum_{j=1}^{[n/2]}\,
\frac{(1-q^{n})\cdots(1-q^{n-j+1})}
{(1-q)\cdots(1-q^{j})}\, M_{n-2j}.\label{overpgen}
\end{align}
These are (continuous) {\em $q$\~Hermite polynomials}
introduced by Szeg\"o and considered in many works; see, e.g.
\cite{ASI}.

Due to (\ref{barpform}), the composition
$\Psi'=q^{-1/4}(1+T)X\pi$ is the raising operator
for the $P$\~polynomials. Namely,
upon its restriction to the symmetric polynomials:  
\begin{align}\label{raisingbar1}
&q^{\frac{n}{2}}\r(P_n)=P_{n+1}\for
\r\equal\frac{X^2\Ga^{-1}-X^{-2}\Ga}
{X-X^{-1}}.
\end{align}
This readily gives (\ref{overpgen}).

The formulas from
(\ref{mubar}), (\ref{constermbar}), (\ref{normppolsbar}) 
read as follows:
\begin{align}\label{normppolsbara1}
&\lan P_{m}(X)P_{n}(X)
\mu_\circ\ran\ 
=\ \de_{mn}\prod_{j=1}^{n}
(1-q^{j})\,,
\end{align}
where $m,n=0,1,\ldots$ and  we set $\mu_\circ=
\mu/\lan\mu_\circ\ran$ for 
\begin{align}
&\mu =  
\prod_{j=0}^\infty (1-X^2 q^{j})
(1-X^{-2}q^{j+1}),\ 
\langle\mu\rangle = 
\prod_{j=1}^{\infty} \frac{1}
{1-q^j}\,.
\label{constermbar1}
\end{align} 
We note that 
$$
\om(\mu)=\mu(Xq^{-1/2};q)=
(-X^2q^{-1})\mu
\and \mu(X^{-1})=-X^{-2}\mu(X).
$$
\smallskip

{\bf Theta functions.}
We set:
\begin{align}\label{gauexp}
\theta\equal &\sum_{n=-\infty}^\infty q^{n^2/4}X^n\\
=&\prod_{j=1}^\infty (1-q^{j/2})(1+q^{\frac{2j-1}{4}}X)
(1+q^{\frac{2j-1}{4}}X^{-1}).\notag
\end{align}
Then $s(\theta)=\theta$ \,and\, 
$\om(\theta)=(q^{-1/4}X)\theta$.\ 
 We will also use 
\begin{align}\label{gauexpu}
\breve{\theta}\equal&\sum_{n=-\infty}^\infty q^{n^2}X^{2n}\\
=&\prod_{j=1}^\infty (1-q^{2j})(1+q^{2j-1}X^2)
(1+q^{2j-1}X^{-2}).\notag
\end{align}
Then 
$s(\breve{\theta})=\breve{\theta}$ and
$\om(\breve{\theta})=(q^{-1}X^2)\breve{\theta}.$\ 
One has (see \cite{ChO})\,:
\begin{align}
&\theta\mu=
\sum_{n=0}^\infty q^{n(n+2)/12}(X^{n+2}-X^{-n}) \,\for
n\neq 2 \hbox{ mod } 3,\notag\\
\label{mugaone}
&\lr P_n P_m\theta\mu\rr= q^{\frac{(m-n)^2}{4}}, \where n,m\ge 0,\\
&\breve{\theta}\mu=
\sum_{n=0}^\infty q^{n(n+1)/3}(X^{2n+2}-X^{-2n}) \for
n\neq 1 \hbox{ mod } 3,\notag\\
\label{mugaoneu}
&\lr P_n P_m\breve{\theta}\mu\rr= 
q^{\frac{(m-n)^2}{4}} \for n-m\in 2\Z \and 0 \hbox{\, otherwise}.
\end{align}

\subsection{Theta-products}\label{sec:Theta-Products}
Our general aim is to
expand the products of $\theta$ and
$\breve{\theta}$ in 
terms of the $P$\~polynomials. We proceed by 
induction using (\ref{mugaone},\ref{mugaoneu}). Generally,
$$
f(X)=\sum_{n=0}^\infty\frac{\lr f\,P_n(X)\mu_\circ\rr}
{\lr P_n^2\mu_\circ\rr}\,P_n \for \lr P_n^2\mu_\circ\rr
\hbox{\ from (\ref{normppolsbara1})}
$$
and for any $s$\~invariant Laurent series $f(X)$.

Let $(\mathbf{A};q)_\infty=\prod_{j=1}^p\prod_{i=0}^\infty
(1-A_jq^i)$ for $\mathbf{A}=\{A_i,\,i=1,\ldots,p\}$.
For instance, $(q;q)_\infty=\prod_{i=1}^\infty
(1-q^i)$. Accordingly,
\begin{align}\label{qn-not}
(\mathbf{A};q)_{\mathbf{n}}\equal\prod_{i=1}^p\prod_{j=0}^{n_i-1}
(1-A_jq^j)=\prod_{i=1}^p (A_i;q)_{n_i},
\end{align}
where  $\mathbf{n}=\{n_i\ge 0,\, 1\le i\le p \}$, $(a;q)_0=1.$
In particular, $(q)_{\mathbf{n}}=
\prod_{i=1}^p\prod_{j=1}^{n_i}(1-q^j)$. 
All expressions below are considered as series in terms
of nonnegative powers of $q$; analytically, one can assume 
that $|q|<1$.
\smallskip

First of all,
\begin{align}\label{thetaform}
&\frac{\th(X)}{(q)_\infty}=\sum_{n=0}^\infty
q^{\frac{n^2}{4}}\frac{P_n}{(q)_n},\\
\label{thetaformu} 
&\frac{\breve{\th}(X)}{(q)_\infty}=\sum_{n=0}^\infty
q^{n^2}\frac{P_{2n}}{(q)_{2n}}.
\end{align}
It will be combined with (\ref{mugaone},\ref{mugaoneu})
in the following particular case of Theorem \ref{MAINTHM}.

Let $\ep(a)=a\mod 2$\, for \,$a\in \Z$,\ 
$\de^\ep_{a,b}=1$ \,if\, $\ep(a)=\ep(b)$ and $0$ otherwise.
For 
$\mathbf{n}=\{n_1,\ldots,n_p\}\subset \Z_+$, we set
$$
\ep^{(\mathbf{n})}_{i,j}=\ep(n_i-n_j)
\for 0\le i,j\le p, \where n_0\equal 0.
$$

\begin{theorem} \label{thetamain}
Let us take a sequence $\boldsymbol{\xi}$
$=\{\xi_i=0,\Join,\ 1\le i \le p\}$,
setting $\theta_i=\theta$ if\, $\xi_i=\Join$\, and \,
$\theta_i=\breve{\theta}$ for $\xi_i=0$. Thus 
$\varpi=\{0\}$ for $\xi=0$ and 
$\varpi=\Join=\Z_2$ for $\xi=\Join$.

For 
$\mathbf{n}=\{n_i\in \Z_+, 1\le i\le p\}$,
\begin{align}\label{pggaga}
&\frac{\prod_{i=1}^p{\theta_i(X)}}{(q)_\infty^p}\ \,=\,\ 
\sum_{\mathbf{n}}\frac{
q^{\bigl(n_1^2+(n_1-n_2)^2+\ldots+(n_{p-1}-n_p)^2\bigr)/4}}
{(q)_\mathbf{n}}
P_{n_p}(X),\\
&\hbox{subject to \,}\ep^{(\mathbf{n})}_{i,i-1}=0 
\hbox{\, if\, } \xi_i=0
\hbox{\, for \,} 1\le i\le p;\ \ep^{(\mathbf{n})}_{1,0}=\ep(n_1). 
\notag
\end{align}
In particular for any fixed $n_p\in \Z_+$,
the corresponding subsum on the right-hand side of (\ref{pggaga}) 
depends only 
on the number of indices $i$ such that $\xi_i=0$ but
not on their specific order in the sequence $\boldsymbol{\xi}$. 
\end{theorem}
\smallskip

{\bf Calculating the coefficients.}
The simplest example of nontrivial combinatorial identities
obtained by permuting $\{\xi_i\}$
is for $p=3$ and $\boldsymbol{\xi}=\{0,\Join,\Join\}$, 
$\boldsymbol{\xi}'=\{\Join,0,\Join\}$. 
Considering the coefficient of $P_k$ for $k=0,1,\ldots$, we 
obtain the following identities:
\begin{align*} 
\sum_{\mathbf{n}}\de^\ep_{n_1,0}\frac{
q^{\bigl(n_1^2+(n_1-n_2)^2+(n_2-k)^2\bigr)/4}}
{(q)_{n_1,n_2}}= 
\sum_{\mathbf{n}}\de^\ep_{n_1,n_2}\frac{
q^{\bigl(n_1^2+(n_1-n_2)^2+(n_2-k)^2\bigr)/4}}
{(q)_{n_1,n_2}}.
\end{align*}
It suffices to assume here that $n_2$ is odd in the both sides
and that $n_1$ is even on the left-hand side, correspondingly, 
odd on the right-hand side. As a matter of fact, this holds for 
any fixed $n_2$, which can be deduced from the Euler 
alternating identity (2.2.1) from \cite{MSZ};
see also \cite{And}, (2.2.6).

\smallskip
Now, taking the constant term of (\ref{pggaga}) and
using (\ref{overpgen}), we arrive at the following
identities:
\begin{align}\label{pggagact}
&\lan\frac{\prod_{i=1}^p{\theta_i(X)}X^m}{(q)_\infty^p}\ran\ \,=\,\
\sum_{\mathbf{n}}\frac{
q^{\bigl(n_1^2+(n_1-n_2)^2+\ldots+(n_{p-1}-n_p)^2\bigr)/4}}
{(q)_{n_p/2+m/2}(q)_{n_p/2-m/2}\prod_{i=1}^{p-1}(q)_{n_i}} \\
&\hbox{ subject to \,} \ep^{(\mathbf{n})}_{i,i-1}=0 
\hbox{\ if\ } \xi_i=0\,
(i\ge 0) \and \ep(n_p-m)=0,\notag
\end{align}
where the indices $k$ in $(q)_k$ are assumed nonnegative. 

When $p=1$, we readily obtain a well-known identity
(see, e.g, \cite{Za}, formula (27)):
\begin{align}\label{pggagid}
&\frac{1}{(q)_\infty}= 
\sum_{k,l\ge 0}^\infty\frac{q^{kl}}
{(q)_k (q)_l} \hbox{\, subject to\, } k=l+m 
\hbox{\, for any\, } m\in \Z.
\end{align}
It is called the Durfee rectangle identity.
Its particular case $m=0$ (corresponding to
$m=0$ in (\ref{pggagact})) is
the Euler identity:
\begin{align}\label{pggagaeu}
&\frac{1}{(q)_\infty}\ =\ 
\sum_{n=0}^\infty\frac{
q^{n^2}}
{(q)_n^2}.
\end{align}
See \cite{MSZ}, formula (2.1.6) and \cite{And}, (2.2.8)
(the Cauchy identity).
\smallskip

\subsection{The case p=2} 
For $\xi_i=\Join\ (i=1,2)$, 
\begin{align}\label{pggagact2}
&\frac{\lan\theta^2\ran}{(q)_\infty^2}\ =\ 
\frac{\sum_{i=-\infty}^\infty q^{i^2/2}}{(q)_\infty^2}=
\frac{\prod_{j=0}^\infty (1+q^{j+1/2})^2}{(q)_\infty}
\\ 
&=\ \sum_{\mathbf{n}}\frac{q^{(n_1^2+(n_1-n_2)^2)/4}}
{(q)_{n_1} \,(q)_{n_2/2}^2},\where \ep(n_2)=0.\notag
\end{align}

For an arbitrary $m\in \Z$,
the coefficient of $X^{m}(m\ge 0)$ here reads:
\begin{align}\label{pggagact2m}
&\frac{\lan X^{m}\theta^2\ran}{(q)_\infty^2}\ \,=\ \, 
q^{m^2/2}\,\frac{\sum_{i=-\infty}^\infty q^{i^2/2}}
{(q)_\infty^2}
\\ 
=\ &\sum_{\mathbf{n}}\de^\ep_{n_2,m}\,
\frac{q^{(n_1^2+(n_1-n_2)^2)/4}}
{(q)_{n_1} \,(q)_{n_2/2+m/2}\,(q)_{n_2/2-m/2}}\,,\notag
\end{align}
assuming that $n_2\pm m\ge 0$.
\smallskip

{\bf Constant term with mu.}
Following Corollary \ref{MAINCOR},
let us now apply $\lan\, \cdot\, \mu_\circ\,\ran$
to (\ref{pggaga}) for $p=2$ and, correspondingly,
for $\boldsymbol{\xi}=\{\Join,\Join\}$ and 
$\boldsymbol{\xi}=\{0,0\}$.
One has: 
\begin{align}\label{pggagact2mu}
&\frac{\lan\theta^2\mu_\circ\ran}{(q)_\infty^2}\, =\,
\sum_{n=0}^\infty\frac{q^{n^2/2}}
{(q)_{n}},\ \,
\frac{\lan\breve{\theta}^2 \mu_\circ\ran}{(q)_\infty^2}\ =
\sum_{n=0}^\infty\frac{q^{2n^2}}
{(q)_{2n}}.
\end{align}
Using (\ref{mugaone})
and formulas (2.8.10) and (2.8.9) from \cite{MSZ},
\begin{align}\label{pggagact2mup}
&\frac{\lan \breve{\theta}^2\mu_\circ\ran}{(q)_\infty^2}\
=\frac{\sum_{j=-\infty}^\infty q^{4j^2-j}}{(q^2\,;\,q^2)_\infty}
=\frac{(-q^3,-q^5,q^8;q^8)_\infty}{(q^2\,;\,q^2)_\infty}\,.
\end{align}

For $\varpi=\{1\}$ (i.e. for $\xi=1$), we will denote the
corresponding $\th_{\varpi}$ by $\hat{\th}$. Then
\begin{align}\label{pggagact2mupp}
&\frac{\lan\hat{\theta}^{\,2}\mu_\circ\ran}{(q)_\infty^2}=
\frac{\lan\theta^2\mu_\circ\ran}{(q)_\infty^2}-
\frac{\lan \breve{\theta}^2\mu_\circ\ran }{(q)_\infty^2}=
q^{1/2}\sum_{n=0}^\infty\frac{q^{2n(n+1)}}{(q)_{2n+1}}\\
&=\,q^{1/2}\frac{\sum_{j=-\infty}^\infty 
q^{4j^2-3j}}{(q^2\,;\,q^2)_\infty}
\ =\ q^{1/2}\frac{(-q,-q^7,q^8;q^8)_\infty}
{(q^2\,;\,q^2)_\infty}\,.\notag
\end{align}

\smallskip

{\bf Modular invariance.}
The relations (\ref{pggagact2mupp},\ref{pggagact2mupp})
are modulo $8$ counterparts of the Rogers-Ramanujan identities,
which are modulo $5$.
They can be found in Table 1 at pg. 44 in \cite{Za};
see also Theorem 3.3 from \cite{VZ}.
There are 3 entries there for $A=1$. 
The first gives that
\begin{align}\label{etazagier}
&\frac{\lan\theta^2\mu_\circ\ran}{(q)_\infty^2}\, =\, 
\sum_{n=0}^\infty\frac{q^{n^2/2}}
{(q)_{n}}\,=\,q^{1/48}\eta(z)^2/(\eta(z/2)\eta(2z))\\
&\hbox{for\ \,} \eta(z)=q^{1/24}\prod_{i=1}^\infty (1-q^n),\where 
q=e^{2\pi \imath z}.\notag
\end{align}

The second entry in this table is directly connected with our
$$
\frac{\lan\theta^2\mu_\circ (X+X^{-1})\ran}{(q)_\infty^2}\, =\, 
q^{1/4}\sum_{n=0}^\infty\frac{q^{(n^2-n)/2}}
{(q)_{n}}\, = \,2q^{1/4-1/24}\eta(2z)/\eta(z).
$$

We see that the function 
$f_{\Join}=q^{-1/48}\frac{\lan\theta^2\mu_\circ\ran}{(q)_\infty^2}$
is modular invariant with respect to $\Ga(2)$
extended by $S$={\tiny
$\left(\begin{array}{cc}
0 & 1 \\
-1 & 0 \\
\end{array}\right)$}.
The action of 
{\tiny $\left(\begin{array}{cc}
a & b \\
c & d \\
\end{array}\right)$}
\, from $SL(2,\Z)$ is $z\mapsto \frac{az+b}{cz+d}$\, and
\begin{align*}
&\Ga(M)=\{A=(a_{ij})\,|\, a_{ij}\in \de_{ij}+M\Z\},
A\in SL(2,\Z),\\
&\Ga_0(M)\ =\ \{\,A\ |\ a_{ij}\in \de_{ij}+M\Z \for 
\{ij\}\neq \{21\}\,\}.
\end{align*}

To be more exact, $f_{\Join}$ is strictly invariant with 
respect to $z\mapsto -1/z$ and 
$f_{\Join}(z+2)=e^{-\frac{\pi \imath}{12}}f_{\Join}(z).$ 
It is in contrast to 
$$
f_0=q^{-1/48}\frac{\lan\breve{\theta}^2
\mu_\circ\ran}{(q)_\infty^2}\,,\ 
f_1=q^{1/2-1/48}\frac{\lan\hat{\theta}^2\mu_\circ\ran}{(q)_\infty^2},
$$
because these two functions are only $\Ga_0(16)$\~invariant
(up to a finite character).
Moreover,
$$
f_0(1/(8z+1)=e^{2\pi \imath/6}f_1(z)\and 
f_1(1/(8z+1)=e^{2\pi \imath/6}f_0(z),
$$
which matches $f_{\Join}(1/(2z+1))=
e^{\frac{\pi \imath}{12}}f_{\Join}(z)$
because $f_{\Join}=f_0+f_1$.

For an arbitrary $p$ and any collections $\boldsymbol{\varpi}$
and corresponding minuscule $c\in P_-$, Corollary \ref{MODULARCOR}
provides the following upper bound for
the modular invariance of the corresponding series
(upon multiplication by a proper fractional powers of $q$); 
it must be modular invariant at least with respect to
the congruence subgroup $\Ga_0(4(p+2))\cap \Ga(2)$.
Thus this estimate is sharp in the case of $p=2$.

\subsection{The case p=3}\label{sec:casep=3}
Here $N=5$.
Let us allow minuscule $c=-\om_r$ in 
Corollary \ref{MAINCOR}, adding 
$P_0=1$ or $P_1=X+X^{-1}$ to the formula for 
$r=0,1$.  We will consider only the {\em atomic} 
sets $\tilde{k}=\{k\}\subset \Z_2$ for $k=0,1$. The
admissible choices for 
$\boldsymbol{\varpi}=\{\om_1,\om_2,\om_3\}$ are
as follows:
\begin{align*}
&r=0\,:\ \, 0\!0\!0=\{\tilde{0},\tilde{0},\tilde{0}\},\
1\!1\!0=\{\tilde{1},\tilde{1},\tilde{0}\},\\
&r=1\,:\ \, 1\!0\!0=\{\tilde{1},\tilde{0},\tilde{0}\},\
1\!1\!1=\{\tilde{1},\tilde{1},\tilde{1}\}.
\end{align*}
We will not distinguish the sequences that 
can be obtained from each other by permutations
since they result in
coinciding  $\Xi_{\boldsymbol{\varpi}}^{3,r}$. As we discussed
above, there are
nontrivial identities reflecting this coincidence.

Recall our main result. Provided that $u+v+w+r=0 \hbox{ mod }(2)$,
\begin{align}\label{pggagact3}
&\Xi_{uvw}^{3,r}=
\frac{\lan \th_u\th_v\th_w P_r\mu_\circ\ran}
{(\prod_{i=1}^\infty(1-q^i))^3}\ =\ 
\sum_{n_1,n_2\ge 0}\frac{q^{(n_1^2-n_1n_2+n_2^2-n_2 r)/2+
\frac{r^2}{4}}}
{(q)_{n_1} \,(q)_{n_2}}\,,\\
&\hbox{subject to \ \,} \ep(n_1)=u,\, \ep(n_2)=u+v 
\,\for \,u,v,w,r \in \Z_2.\notag
\end{align}
When $\boldsymbol{\varpi}=\Join$,
these formulas are from the last entry of Table 1 in \cite{Za}
for the matrix {\tiny 
$\left(\begin{array}{cc}
1 & -\frac{1}{2} \\
-\frac{1}{2} & 1 \\
\end{array}\right)$}. See also the table from
Theorem 3.4 from \cite{VZ} (the first entry there).
The latter table provides the eta-type formulas for 
$\Xi_{\Join}^{3,0}$ and $\Xi_{\Join}^{3,1}$, but we
need them here for atomic $\varpi$. They are
certain (not quite trivial) splits of the formulas given in \cite{VZ}: 
\begin{align}\label{r-r501}
q^{-1/20}\,\Xi_{0\!0\!0}^{3,0}\,=\,
&(\,\th_{5,\frac{3}{4}}(2z)\,+\,
\th_{5,\frac{13}{4}}(2z)\,)\eta(z)\,/\,
(\eta(2z)\eta(z/2))\\
-&\ \,\th_{5,2}(2z)\eta(2z)/\eta(z)^2,\notag\\
\label{r-r51}
q^{-1/20}\,\Xi_{1\!1\!0}^{3,0}\,=\, 
&\ \,\th_{5,2}(2z)\eta(2z)/\eta(z)^2,\where\\
\th_{5,m}(z)\!\equal\!
\sum_{n}&(-1)^{[n/10]}q^{n^{2}/40}
\for n\in 2m-1+10\,\Z\,,\label{theta5}
\end{align}
\begin{align}
\label{r-r502}
q^{-4/20}\,\Xi_{1\!1\!1}^{3,1}\,=\,
&\th_{5,\frac{3}{2}}(z)\th_{5,2}(2z)\eta(z)^3/
(\eta(z/2)^2\eta(2z)^2)\eta(10z))\\
-&\th_{5,1}(2z)\eta(2z)/\eta(z)^2,\notag\\
\label{r-r52}
q^{-4/20}\,\Xi_{1\!0\!0}^{3,1}\,=\,
&\th_{5,1}(2z)\eta(2z)/\eta(z)^2.
\end{align}
\medskip

{\bf Level-rank duality}.
The summations in the formulas (\ref{pggagact3})
naturally appear in the theory of $\hat{\mathfrak{sl}}_3$.
Namely, our summations divided by $\eta(z)^2$ and
with proper fractional powers of $q$ coincide
with formulas (5.16-19) in \cite{Geo} for certain
level-two {\em string functions}. As a matter of fact,
the string functions listed there are all independent ones
for such a level. We note that $q^{-2/15}$
and $q^{-1/30}$ in \cite{Geo} and below are our 
$q$\~power corrections (necessary to ensure the modular invariance
of $\Xi$ and the string functions)
adjusted due to the division by $\eta(z)^2$\,:\ \,
$$
-\frac{1}{20}=-\frac{2}{15}+\frac{1}{12},\ 
-\frac{4}{20}=-\frac{1}{30}+\frac{1}{12}-\frac{1}{4}.
$$

Using the notations from \cite{Geo}
(see also \cite{KP}),
\begin{align*}
\frac{q^{-2/15+1/12}}
{\eta(z)^2}\,\Xi_{0\!0\!0}^{3,0}\,&=\, c_0^{2\hat{\La}_0},\ 
&\frac{q^{-1/20-1/12}}
{\eta(z)^2}\,\Xi_{1\!1\!0}^{3,0}\,=\, c^{2\hat{\La}_0}_{\al_1+\al_2},
\\ 
\frac{q^{-1/30+1/12-1/4}}
{\eta(z)^2}\,\Xi_{1\!0\!0}^{3,1}\,&=\, 
c_{\La_1}^{\hat{\La}_0+\hat{\La}_1},\ 
&\frac{q^{-1/30+1/12-1/4}}
{\eta(z)^2}\,\Xi_{1\!1\!1}^{3,0}\,=\, 
c_{\La_1+\al_2}^{\hat{\La}_0+\hat{\La}_1}.
\end{align*}

We hope to discuss this duality for arbitrary 
(atomic) $\boldsymbol{\varpi}$ in further works.
\medskip

{\bf Modular invariance.}
Let us briefly discuss the modular properties of these
functions. It is directly connected with the
action of $PSL(2,\Z)$ in the
$4$-dimensional Verlinde algebra for $A_1$ of level $3$;
the corresponding products of $\th_{\boldsymbol{\varpi}}$
form a basis in this space.
This guarantees that the
whole $PSL(2,\Z)$ acts in the linear space 
generated by $\Xi_{uvw}^{3,r}$ from (\ref{pggagact3}). 
Explicitly, this action is as follows.
\smallskip

All four functions are modular invariant 
 with respect to $(T_-)^M$ where $M=10$
(exactly) and under the action of $T_+$
(up to proportionality). The total modular invariance 
subgroup is $\Ga_0(10)$ (up to a character). 

As for $T_-^5$, setting 
$$
{}^\sharp\!\overrightarrow{\Xi}^0\ =\ q^{-1/20}
(\,\Xi_{0\!0\!0}^{3,0},\,\Xi_{1\!1\!0}^{3,0}\,)^{tr},\ \
{}^\sharp\!\overrightarrow{\Xi}^1\ =\ q^{-4/20}
(\,\Xi_{1\!1\!1}^{3,1},\,\Xi_{1\!0\!0}^{3,1}\,)^{tr},
$$
\begin{align}\label{r-r-0}
&(T_-^5)\,\bigl(\,{}^\sharp\!\overrightarrow{\Xi}^0,\,
{}^\sharp\!\overrightarrow{\Xi}^1\,\bigr)\ =\ 
\left(
  \begin{array}{cc}
    -\frac{1}{2} &  \ \,\frac{3}{2} \\
    \ \frac{1}{2} &   \frac{1}{2} \\
  \end{array}
\right)\,
\bigl(\,{}^\sharp\!\overrightarrow{\Xi}^0,\,
{}^\sharp\!\overrightarrow{\Xi}^1\,\bigr).
\end{align}

The eigenvalues
of $T_+$ are correspondingly
$$
\pm e^{2\pi \imath/20}\, \for\,  {}^\sharp\Xi_{1\!0\!0}^{3,1},\  
{}^\sharp\Xi_{1\!1\!1}^{3,1}\ 
\and \ \pm e^{- 2\pi \imath/20}\, \for\, 
{}^\sharp\Xi_{0\!0\!0}^{3,0},\ {}^\sharp\Xi_{1\!1\!0}^{3,0}.
$$
This is actually obvious, since the $q$\~series for 
$q^{-1/4}\Xi_{1\!0\!0}^{3,1}$, $q^{1/4}\Xi_{1\!1\!1}^{3,1}$ and 
$\Xi_{0\!0\!0}^{3,0}$, $q^{-1/2}\Xi_{1\!1\!0}^{3,0}$
contain only integral powers of $q$, and $T_+(q^{1/2})=-q^{1/2}.$
The corresponding eigenvalues will be only due to
the fractional $q$\~powers in (\ref{r-r501},\ref{r-r51}) and
(\ref{r-r502},\ref{r-r52}), i.e. due to the 
$q$\~normalization.
\smallskip

Furthermore, $T_-^2$ preserves the two-dimensional 
spaces
$$
\C\ {}^\sharp\Xi_{1\!0\!0}^{3,1}\,\oplus\,
\C\ {}^\sharp\Xi_{1\!1\!0}^{3,0}\, \and\,
\C\ {}^\sharp\Xi_{0\!0\!0}^{3,0}\,\oplus\,
\C\ {}^\sharp\Xi_{1\!1\!1}^{3,1}.
$$
Setting now ${}^\sharp\overline{\overline{\Xi}}$=
$\bigl(\,{}^\sharp\!\overleftarrow{\Xi}^0,\,
{}^\sharp\!\overleftarrow{\Xi}^0\bigr)$ for
$$
\ \
{}^\sharp\!\overleftarrow{\Xi}^0\ =\ 
(\,q^{-4/20}\Xi_{1\!0\!0}^{3,1},\,q^{-1/20}\Xi_{1\!1\!0}^{3,0}\,)^{tr},
\ \ 
{}^\sharp\!\overleftarrow{\Xi}^1\ =\ 
(\,q^{-4/20}\Xi_{1\!1\!1}^{3,1},\,q^{-1/20}\Xi_{0\!0\!0}^{3,0}\,)^{tr},
$$
\begin{align}\label{r-r-01}
&(T_-^2)\,{}^\sharp\overline{\overline{\Xi}}\,=\, 
\left(\begin{array}{cc}
    \cos(\frac{\pi}{5})\!+\! \imath\frac{\sin(\pi/5)}{\sqrt{5}}& 
                2\imath\frac{\sin(\pi/5)}{\sqrt{5}} \\
                2\imath\frac{\sin(\pi/5)}{\sqrt{5}}&   
    \cos(\frac{\pi}{5})\!-\!\imath\frac{\sin(\pi/5)}{\sqrt{5}} \\
\end{array}\right)\,
{}^\sharp\overline{\overline{\Xi}}.
\end{align}

Since the modular transformations
$T_+$, $T_-^5$ and $T_-^2$ act in this $4$\~dimensional space, 
the whole $PSL(2,\Z)$ acts here, and in a sufficiently explicit
way.

Some of these facts are well known. The functions
$\Xi_{1\!0\!0}^{3,1}$ and $\Xi_{1\!1\!0}^{3,0}$
are directly related to the Rogers-Ramanujan identities
with $q^2$ instead of $q$ and a simple common eta-type factor:
\begin{align}\label{Rog-Ram1}
&
q^{-\frac{1}{4}}\Xi_{1\!0\!0}^{3,1}\!=\!
q^{-\frac{1}{20}}\th_{5,1}(2z)
\frac{\eta(2z)}{\eta(z)^2}\!=\!
\!\!\sum_{n=0}^\infty \frac{q^{2n^2}}{\prod_{j=1}^n(1-q^{2j})}
\prod_{j=1}^{\infty}(1+q^j)^2,\\
\label{Rog-Ram2}
&q^{-\frac{1}{2}}\Xi_{1\!1\!0}^{3,0}\!=\!
q^{-\frac{9}{20}}\th_{5,2}(2z)
\frac{\eta(2z)}{\eta(z)^2}\!=\!
\!\!\sum_{n=0}^\infty \frac{q^{2n^2+2n}}{\prod_{j=1}^n(1-q^{2j})}
\prod_{j=1}^{\infty}(1+q^j)^2.
\end{align}
Compare with the classical formulas:
\begin{align}\label{RRclass}
&G(z)=\!\sum_{n\ge 0}\!\frac{q^{n^2}}{(q)_n}\!=q^{\frac{1}{60}}\,
\frac{\th_{5,1}(z)}{\eta(z)},
H(z)=\!\sum_{n\ge 0}\!\frac{q^{n^2+n}}{(q)_n}\!=q^{-\frac{11}{60}}\,
\frac{\th_{5,2}(z)}{\eta(z)}.
\end{align}
See, e.g. formula (23) from \cite{Za}. 
For instance, combining (\ref{r-r-0}) and
(\ref{r-r-01}) we arrive at formula (24) from \cite{Za},
describing the action of $z\mapsto -1/z$ on $G,H$.

Let us comment on the powers of
$q$ that appear here. The functions $q^{-1/40}\th_{5,1}(z)$
and $q^{-9/40}\th_{5,2}(z)$ are power series in terms of 
integral powers of $q$,
which is obvious from the definition in (\ref{theta5})
and can be immediately seen from (\ref{RRclass}); 
use that $q^{-1/24}\eta(z)$ is such a series. Upon $z\mapsto 2z$,
we readily arrive at the fractional $q$\~powers in 
(\ref{Rog-Ram1},\ref{Rog-Ram2}).
\smallskip

\subsection{A coset interpretation}
We will focus only on formula
(\ref{Rog-Ram2}).
According to Section \ref{Tcosets},
we need to consider three level-one integrable modules 
$\,M_1=L_{1,0},\, M_2=L_{1,0},\, M_3=L_{0,1}\,$ of 
$\hat{\mathfrak{g}}=\hat{\mathfrak{sl}}_2$
with the highest weights $\{1,0\}$, $\{1,0\}$ and $\{0,1\}$
in the standard notation. The module $M_3=L_{0,1}$
is the so-called vacuum representation, where 
$\mathfrak{sl}_2\otimes\C[t]$
acts as zero. Then
$\nu(L\,;\, M_1,M_2,M_3)= 
\hbox{Hom\,}_{\hat{\mathfrak{sl}}_2}\,
(L\,, M_1\otimes M_2\otimes M_3)$ is a natural
representation of the coset algebra defined for the diagonal
embedding 
$\hat{\mathfrak{sl}}_2\hookrightarrow 
\hat{\mathfrak{sl}}_2\times\hat{\mathfrak{sl}}_2\times
\hat{\mathfrak{sl}}_2$. 

Here we 
consider $\,M_1\otimes M_2\otimes M_3\,$ as a submodule of
$\,M^{\otimes 3}\,$ for $\,M=L_{1,0}\oplus L_{0,1}$. The passage
to the adjoint (graded) module followed by taking the constant
term with $\mu$ is standard here. It is more subtle for the 
identity from (\ref{Rog-Ram1}), where we need to multiply
the corresponding character by $P_{1}$ before taking the
constant term.

The level-rank duality says that the coset 
$(\hat{\mathfrak{sl}}_2\times\hat{\mathfrak{sl}}_2\times
\hat{\mathfrak{sl}}_2)/\hat{\mathfrak{sl}}_2$ 
(all are of level $1$) is the same 
as the coset of $\hat{\mathfrak{sl}}_3$ at the level $2$ over
the Heisenberg algebra 
$$
\hat{\mathfrak{h}}_3=\mathfrak{h}_3\otimes \C[t,t^{-1}]
\hbox{\ \,for the Cartan subalgebra \ \,} 
\mathfrak{h}_3\subset \mathfrak{sl}_3.
$$ 
The coset 
$(\hat{\mathfrak{sl}}_3/\hat{\mathfrak{h}}_3)$ can be naturally
considered as a product of two level-$2$ cosets 
$(\hat{\mathfrak{sl}}_3/\hat{\mathfrak{gl}}_2)$ and
$(\hat{\mathfrak{sl}}_2/\hat{\mathfrak{h}}_2)$, where
we denote by 
$\hat{\mathfrak{h}}_2$ the Heisenberg algebra in 
$\hat{\mathfrak{sl}}_2$ and use the standard embedding
$\mathfrak{gl}_2 \hookrightarrow\mathfrak{sl}_3$.
This presentation is simply because of consecutive 
taking the invariants.
Since the level is $2$, these two cosets belong correspondingly
to the Virasoro $\{4,5\}$ and $\{3,4\}$ minimal models.
\smallskip

The term $\prod_{j=1}^{\infty}(1+q^j)$ 
on the right-hand side of (\ref{Rog-Ram2})
corresponds to the Virasoro module of weight $1/16$ from
the $\{3,4\}$ minimal model. In the standard
notation, it is $\varphi_{1,2}$.

The remaining term 
$$
\prod_{j=1}^{\infty}(1+q^j)\,
\sum_{n=0}^\infty q^{2n^2+2n}/(\prod_{j=1}^n(1-q^{2j}))
$$ 
from (\ref{Rog-Ram2}) is of more involved nature. Namely, it is 
the character 
of the Virasoro module $\varphi_{2,2}$ from the minimal model
$\{4,5\}$.

For this identification, we use that the Virasoro
$\{4,5\}$ minimal model is related to the minimal models for the
Neveu-Schwarz and Ramond superalgebras. For instance,
it gives that the Virasoro module $\varphi_{2,2}$ can be 
viewed as an irreducible representation of the Ramond algebra.
The latter algebra has a basis $\{L_i,S_j\}$, where $L_i$ are
Virasoro generators and $S_j$ are odd satisfying
$[S_i,S_j]_+=L_{i+j}$. Thus knowing the action of $\{S_j\}$
is sufficient. One can check that  $\varphi_{2,2}$
has a monomial basis 
\begin{align*}
&\{\,S_{-1}^{a_1}S_{-2}^{a_2}\cdots S_{-m}^{a_m}\,\},
\where 0\le a_i<4 \and \\
&\hbox{if\ \,} a_i>1,
\hbox{\  then\,  both\ \,} a_{i-1}\le 1 \hbox{ and }  a_{i+1}\le 1.
\end{align*}
This results in the required formula for the
character and gives a coset interpretation of
(\ref{Rog-Ram2}). The remaining three identities 
from (\ref{r-r501},\ref{r-r502},\ref{r-r52}) can be also 
interpreted in this way, which is quite interesting in 
coset theory.
\smallskip

\subsection{Theta of level 1/2}\label{sec:half-int}
According to our analysis, the ``main stream" of the theory
of Rogers-Ramanujan type identities (for instance, papers
\cite{An,War}) can be connected with ``square roots"
of level-one theta functions. Stimulated by formulas
(\ref{warnaarB}) and (\ref{warnaarC}), we
defined and calculated (numerically) the $L$\~sums upon 
division of the  $A$\~matrices considered above by $2$.
They appeared rational for all root systems (see below).
We will denote them by $L_R^\flat$. 


It is explained in \cite{Nak2} that these formulas can be
obtained from known identities using formula (1.7) from \cite{Nak} 
for $\ell=3$; see Corollary 1.9 there.
In the $BCFGT$\~cases, the so-called folding construction 
can be used; see Section 9 from \cite{IKNS} and \cite{Nak2}. 
Importantly, the $Y$\~systems establish a {\em direct} link 
to the Langlands functoriality of the {\em affine} root systems.
We will touch it upon below, when discussing functoriality
properties of the effective central charges.
\smallskip 

Let $c_{\hbox{\tiny eff}}=L_R^\flat/h\,$ for 
the Coxeter number $h$; notice that it was $2 L/h$ before.
The matrix $A_R^\flat= A_R/2=(a_{ij}^\flat)$ 
has the entries $a_{ij}^\flat=(\om_i,\om_j)$. Recall that the weight 
lattice $P$ is supplied with the standard form 
$(\cdot,\cdot)$ normalized
by the condition $(\al_{\sht},\al_{\sht})=2$. 

In terms of
$\nu_i=\frac{(\al_i,\al_i)}{2}$, the $Q$\~system is
as follows\,:  
\begin{align}\label{Qsystem}
&(1-Q_i)^{\nu_i}
=\prod_{j=1}^n Q_j^{a'_{ij}}\, (1\le i\le n),\ 
L'_R = \frac{6}{\pi^2}\,\sum_{i=1}^n \nu_i
L(Q_i)\,;  
\end{align}
$A'$ here is $A$ or $A^\flat$.  
The $Q$\~system in the case of $T_n$ identically coincides 
with that for $A_{2n}$ (with ${}^\flat$ or without) 
with a reservation that the
number of terms in $L'_{T_n}$ is $n$ versus $2n$ for $L'_{A_{2n}}$.
The $T$\~type $Q$\~system is defined for $C_n$, but without using 
$\nu_i$ in (\ref{Qsystem}).
\smallskip

We note that a somewhat different
$Q$\~system naturally appears when applying the method from \cite{VZ}
(and previous works):
\begin{align}\label{Qsystemd}
&1-Q_i
=\prod_{j=1}^n Q_j^{a'_{ij}/\nu_j}\, (1\le i\le n),\ 
\tilde{L}'_R \,=\, \frac{6}{\pi^2}\,
\sum_{i=1}^n \frac{\nu_{\lng}}{\nu_i}
L(Q_i)\,,
\end{align}
where $A'$ is $A$ or $A^\flat$.  It is simple to see that 
$\tilde{L}'_R=L'_{R^\vee}$ for $R=B_n,\,C_n,\,F_4,\,G_2.$
\smallskip

{\bf The table.}
Let us compare the values of $L$ and the
corresponding $c_{\hbox{\tiny eff}}$ without
and with $\flat$. To avoid misunderstanding with the 
normalization, let us list the coefficients $a_{11}^\flat$.
They are $n/(n+1)$ for $A_n$, $4/3$ for $E_6$,
$2$ for $B_n,E_7,G_2$ and $1$ otherwise.
The values of $L$ and effective central charges are
as follows:
\medskip

{\footnotesize
\begin{tabular}{|c|c|c|c|c|c|}
  \hline
  $R_n$ & $L_R$ & $c_{\hbox{\tiny eff}}$ && $L_R^\flat$ & 
$c_{\hbox{\tiny eff}}^\flat$ \\   
\hline
  $A_n$ & n(n+1)/(n+3)  & 2n/(n+3)     &&  
n(n+1)/(n+4)     & n/(n+4)   \\   
  $B_n$ & n(2n-1)/(n+1) & (2n-1)/(n+1) && 
2n(2n-1)/(2n+3)  & (2n-1)/(2n+3)   \\   
  $C_n$ & n     & 1    && 2n(n+1)/(2n+3)  & (n+1)/(2n+3) \\   
  $D_n$ & n-1   & 1    && 2(n-1)n/(2n+1)  & n/(2n+1)     \\   
  $E_6$ & 36/7  & 6/7  && 24/5  & 2/5  \\   
  $E_7$ & 63/10 & 7/10 && 6     & 1/3  \\   
  $E_8$ & 15/2  & 1/2  && 80/11 & 8/33 \\   
  $F_4$ & 36/7  & 6/7  && 24/5  & 2/5  \\   
  $G_2$ & 3     &  1   && 8/3   & 4/9  \\   
  $T_n$ & n(2n+1)/(2n+3)  &  2n/(2n+3) && 
n(2n+1)/(2n+4)   & n/(2n+4)            \\         
\hline
\end{tabular}
}   
  
\smallskip
\centerline {\em Table 1. Values of $L_R, L_R^\flat$ 
and $c_{\hbox{\tiny eff}}$, $c_{\hbox{\tiny eff}}^\flat$}
\medskip

See \cite{KM} about the history, physics meaning, merits and demerits
of $c_{\hbox{\tiny eff}}$ in the $A$-$D$-$E$ elastic scattering 
theories. Presumably it measures the
massless degrees of freedom of a theory, somehow analogous to the 
calculation we peformed above in several examples
of the ``fraction of the missing terms" in 
the numerator of product formulas. One of the 
problems is that $c_{\hbox{\tiny eff}}$ (without $\flat$)
do not reflect well the connections between different root systems. 
Also, the effective central charges for the series $D_n$ are not 
what one could expect. This problem seems to be better
addressed when using $A^\flat=A/2$. 

Indeed, the coincidences and other relations of the effective 
central charges in the $\flat$\~case  become practically
always meaningful in Kac-Moody theory and in the theory of
coset models. 
The effective charges are equal to each other for the 
following pairs:
$$
D_{n+1}\leftrightarrow C_n,\  B_n \leftrightarrow A_{2n-1},\ 
E_6 \leftrightarrow F_4 \leftrightarrow D_2,\ 
G_2\leftrightarrow D_4,\  
E_7 \leftrightarrow A_2.
$$  
The relation $\,2(c_{\hbox{\tiny eff}}^\flat)_{E_8}=
(c_{\hbox{\tiny eff}}^\flat)_{D_{16}}$ also makes sense
from the viewpoint of representation theory. 
Almost all of these relations are directly related to the
symmetries of the $Q$\~systems (and the uniqueness of the
solutions in the range $\{0<Q_i<1\}$) and to the Langlands 
duality for the affine root systems. See \cite{Nak2} and references
therein. Moreover, let us mention that
the $A_n$ central charges are generally greater than $1/2$ and
intersect the $D$\~charges, strictly smaller
than $1/2$, only at $D_1$ and $D_3$. A similar relation
holds for $B$ and $C$. 
\smallskip

For the physics part, the breakthrough paper \cite{Zam}
(the three-state Potts model, $c=\frac{4}{5}$) and its further 
developments for the 
$W(A_n)$\~algebras result in  $c=2n/(n+3)$ (see \cite{KM}).
It is the central charge of the unitary CFT
describing the $\Z_{n+1}$\~parafermions, which coincides
with $(c_{\hbox{\tiny eff}})_{A_n}$. In unitary theories,
$d_0=0$ so $c_{\hbox{\tiny eff}}=c$.
In the $\flat$\~case, the central charge for $A_n$
is $n/(n+4)$ (maybe up to a common coefficient of proportionality).
\smallskip

Our approach, based on expanding products of 
level-one starting theta functions in terms of the $q$\~Hermite
polynomials, cannot be directly extended to the $\flat$\~case.
Theta functions of level $1/2$ and the corresponding
Kac-Moody modules are needed here. The latter modules are not available
in integrable Kac-Moody theory, but the Virasoro and
$W$\~algebras, parafermions, Neveu-Schwarz and Ramond 
superalgebras provide certain substitutes. DAHA seems to have
greater flexibility here, though this cannot be translated
to Kac-Moody theory (so far).
\smallskip

{\bf Nakanishi's note}.
Let us comment on \cite{Nak2}; we thank Tomoki Nakanishi for 
various discussions. It suggests that one can try to obtain
the $\flat$ case by considering $p=3$ 
(coinciding with $\ell$ in \cite{Nak2})
combined with the folding construction. Indeed, 
following his note, the approach via Nahm's conjecture gives the 
rationality of $\tilde{L}^\flat_R$. However, one needs  
the theory of $Y$\~systems to obtain these rational numbers.
It is challenging that such a reduction (from $2n$ variables
for $p=3$ to $n$ variables using the symmetry of $A_2$) cannot be 
seen at the level of the corresponding Rogers-Ramanujan identities. 

Let us provide some details. For an arbitrary root system $R$,
taking $\boldsymbol{\varpi}=\{\Join,\,\Join,\,\Join\}$ and 
$c=0$ for $p=3$,
the series from (\ref{pggmcor}) reads as follows:
\begin{align}\label{pggmcor3}
&
\sum_{b_1,\ldots,b_2}
\frac{q^{\bigl(b_1^2+(b_1-b_2)^2+(b_{2})^2)/2}}
{\prod_{i=1}^{2}\prod_{j=1}^n\prod_{k=1}^{-(\al_j^\vee,b_i)}
(1-q_j^{k})}\,,
\end{align} 
where the summation is over all possible $b_1,b_2\in P_-$.
If one formally makes here $b_1=b=b_2$, then it becomes
\begin{align}\label{pggmcor33}
&\sum_{b}
\frac{q^{b^2}}
{\prod_{j=1}^n\prod_{k=1}^{-(\al_j^\vee,b)}
(1-q_j^{k})^2}\,,
\end{align} 
i.e. with the same $q$\~powers in the numerators 
and the squares of the denominators versus the 
non-$\flat$ case with $p=2$.

We do not see direct relations between the modular properties
of (\ref{pggmcor33}) and (\ref{pggmcor3}),
but one can proceed as follows. The uniqueness Lemma 2.1 
from \cite{VZ} implies that the canonical solution of the 
$Q$\~system for (\ref{pggmcor3}) in the range $\{0<Q_i<1\}$ is
invariant under the transposition $b_1\leftrightarrow b_2$. 
Therefore it also solves the $Q$\~systems for (\ref{pggmcor33}), 
which is the $\flat$\~ones due to the squares in the 
denominators of (\ref{pggmcor33}). 

More generally, the $Q$\~systems can be reduced from 
$(p-1)n$ variables to $[p/2]n$ variables 
for any given root system $R$ and level $p$;
we impose the symmetry corresponding to the standard 
automorphism of $A_{p-1}$.  
Let as assume that Nahm's conjecture holds in this
situation and that there is a similar reduction for
all (complex) solutions of the $Q$\~system; the
uniqueness claim guarantees it only in the range $\{0<Q_i<1\}$.
Then this makes the modular invariance of the Rogers-Ramanujan 
series in the $\flat$\~case for $p-1$ a (conditional) corollary
of that in the non-$\flat$\~case for $p$. We hope to extend the DAHA 
approach to obtain this fact without using Nahm's conjecture.
\smallskip

In conclusion, we note that
the exact rational values of the $L$\~sums can be generally
obtained from the product formulas or similar identities for the 
modular Rogers-Ramanujan series. DAHA methods can be 
used here, but we touch the product-type formulas only a little in 
this paper. Therefore such identities can be viewed as 
``quantization" of the $Q$\~systems and the corresponding 
dilogarithm formulas. Nahm's conjecture outlines the class of  
$Q$\~systems that can be ``quantized," i.e. lifted to modular 
invariant Rogers-Ramanujan type series. Continuing 
this line, one can try to interpret the coordinate Bethe ansatz,
associated with unitary DAHA modules and analytic properties
of the eigenfunctions of the corresponding QMBP,  
as some kind of quantization of TBA and the 
$Y$\~systems, but this is beyond the present paper.

\vskip -0.4cm
\medskip
\bibliographystyle{unsrt}

\end{document}

%% file: macro.tex
\renewcommand{\tilde}{\widetilde}
\renewcommand{\hat}{\widehat}

\newcommand{\Z}{{\mathbb Z}}
\newcommand{\Q}{{\mathbb Q}}

\newcommand{\C}{{\mathbb C}}
\newcommand{\R}{{\mathbb R}}

\def\HH{\mbox{${\mathcal H}$\kern-5.2pt${\mathcal H}$}}

\newtheorem{theorem}{Theorem}[section]

\newtheorem{proposition}[theorem]{Proposition}
\newtheorem{definition}[theorem]{Definition}
\newtheorem{lemma}[theorem]{Lemma}
\newtheorem{corollary}[theorem]{Corollary}

\newtheorem{theorem }{Theorem}[section]
\newtheorem{maintheorem }[theorem]{Main Theorem}
\newtheorem{proposition }[theorem]{Proposition}
\newtheorem{definition }[theorem]{Definition}
\newtheorem{lemma }[theorem]{Lemma}
\newtheorem{corollary }[theorem]{Corollary}
\newtheorem{notation }[theorem]{Notation}
\newtheorem{remark }[theorem]{Remark}
\newtheorem{example }[theorem]{Example}

\newtheorem{ maintheorem }[theorem]{Main Theorem}
\newtheorem{ theorem}{Theorem}[section]
\newtheorem{ proposition}[theorem]{Proposition}
\newtheorem{ definition}[theorem]{Definition}
\newtheorem{ lemma}[theorem]{Lemma}
\newtheorem{ corollary}[theorem]{Corollary}
\newtheorem{ notation}[theorem]{Notation}
\newtheorem{ remark}[theorem]{Remark}
\newtheorem{ example}[theorem]{Example}

\def\for{\  \hbox{ for } \ }
\def\iif{ \ \hbox{ if } \ }
\def\when{ \ \hbox{ when } \ }
\def\where{\  \hbox{ where } \ }
\def\and{\  \hbox{ and } \ }
\def\and{\  \hbox{ and } \ }

\def\equal{\stackrel{\,\mathbf{def}}{= \kern-3pt =}}

\def\la{\lambda}
\def\La{\Lambda}
\def\om{\omega}

\def\th{\theta}
\def\al{\alpha}

\def\ga{\gamma}
\def\ep{\epsilon}

\def\de{\delta}

\def\si{\sigma}

\def\Ga{\Gamma}
\def\ze{\zeta}

\def\vep{\varepsilon}

\def\vth{{\vartheta}}

\def\tal{\tilde{\alpha}}

\def\tGa{\tilde{\Gamma}}

\def\tw{\widetilde w}
\def\tW{\widetilde W}

\def\tz{\tilde z}
\def\tb{\tilde b}

\def\tR{\tilde R}

\def\hw{\widehat{w}}
\def\hW{\widehat{W}}

\def\hv{\hat{v}}

\def\0{\mathbf{0}}

\def\f{\mathcal{F}}
\def\çF{\mathcal{F}}

\def\r{\mathcal{R}}
\def\l{\mathcal{L}}

\def\a{\mathcal{A}}
\def\h{\mathcal{H}}

\def\v{\mathcal{V}}

\def\w{\mathcal{W}}

\def\b{\mathcal{B}}

\def\lan{\langle}

\def\ran{\rangle}

\def\lng{\hbox{\rm{\tiny lng}}}
\def\sht{\hbox{\rm{\tiny sht}}}

\newcommand{\sq}{\phantom{1}\hfill$\qed$}

\newcommand{\lr}{\langle}
\newcommand{\rr}{\rangle}

\newcommand{\sgn}{\mbox{sgn}}

\def\HH{\mathfrak{H}}

\def\HH{\hbox{${\mathcal H}$\kern-5.2pt${\mathcal H}$}}

\font\smm=msbm10 at 12pt 
\def\symbol#1{\hbox{\smm #1}}
\def\lsmash{{\symbol n}}

\def\#{\sharp}
